\documentclass[a4paper,leqno,12pt]{amsart}
\usepackage{amsfonts}

\usepackage{graphicx}
\usepackage{amsmath}

\newtheorem{theorem}{Theorem}

\newtheorem{corollary}[theorem]{Corollary}

\newtheorem{definition}{Definition}

\newtheorem{lemma}[theorem]{Lemma}

\newtheorem{proposition}[theorem]{Proposition}
\newtheorem{remark}{Remark}

\begin{document}
\title[Generalized Gevrey ultradistributions]{Colombeau generalized Gevrey ultradistributions and their microlocal analysis}
\author{Khaled BENMERIEM}
\author{Chikh BOUZAR}
\address{Centre Universitaire de Mascara. Algeria }
\email{benmeriemkhaled@yahoo.com }
\address{Department of Mathematics, Oran-Essenia University. Algeria }
\email{bouzar@univ-oran.dz ; bouzar@yahoo.com}
\date{}
\maketitle

\begin{abstract}
The purpose of this paper is to construct and to study algebras of
generalized Gevrey ultardistributions. We define the generalized Gevrey wave
front and give its main properties. As a fundamental application, the well
known H\"{o}rmander's theorem on the product of two distributions is
extended to the case of generalized Gevrey ultradistributions
\end{abstract}

\section{Introduction}

The nonlinear theory of generalized functions initiated by J. F. Colombeau,
\cite{Col1} and \cite{Col2}, in connection with the problem of
multiplication of Schwartz distributions \cite{Schw}, has been developed and
applied in nonlinear and linear problems, \cite{Biag}, \cite{Col3}, \cite
{Ober} and \cite{NPS}. The recent book \cite{GKOS} gives further
developments and applications of Colombeau generalized functions. Some
methods of constructing algebras of generalized functions of Colombeau type
are given in \cite{A-R} \ and \cite{Mart}. The proceedings \cite{DHMV} and
\cite{GHKO} present different results on nonlinear analysis of Colombeau
generalized functions.

Ultradistributions, important in theoretical as well applied fields, see
\cite{Kom}, \cite{LM} and \cite{Rod}, are natural generalization of Schwartz
distributions, so it is natural to search for algebras of generalized
functions containing ultradistributions, to study and to apply them. This is
the purpose of this paper.

We first introduce algebras of generalized Gevrey ultradistributions, such a
question is considered in the only papers \cite{Gram}, \cite{Pilip} and \cite
{DHPV}. We then develop a Gevrey microlocal analysis suitable for these
algebras in the spirit of \cite{Hor}, \cite{Rod} and \cite{NPS}. Finally, we
give an application through a generalization of H\"{o}rmander's theorem on
the wave front of the product of two distributions, this is also an
extension of the result of \cite{Hor-Kun}.

The algebras $\mathcal{G}^{s}\left( \Omega \right) $ of generalized Gevrey
ultardistributions are represented by nets of smooth\ functions \ $%
f_{\varepsilon }$\ with exponential growth\ in $\varepsilon $ depending on
the Gevrey order $s$, more precisely
\begin{equation*}
\mathcal{G}^{s}\left( \Omega \right) =\frac{\mathcal{E}_{m}^{s}\left( \Omega
\right) }{\mathcal{N}^{s}\left( \Omega \right) }\text{ \ ,}
\end{equation*}
where $\mathcal{E}_{m}^{s}(\Omega )$ is the space of $\left( f_{\varepsilon
}\right) _{\varepsilon }\in C^{\infty }\left( \Omega \right) ^{\left] 0,1%
\right[ }$ satisfying for every compact subset $K$ of $\Omega $, $\forall
\alpha \in \mathbb{Z}_{+}^{n},\exists C>0,\exists k>0,$ $\exists \varepsilon
_{0}\in \left] 0,1\right[ ,$
\begin{equation}
\left| \partial ^{\alpha }f_{\varepsilon }\left( x\right) \right| \leq C\exp
\left( k\varepsilon ^{-\frac{1}{2s-1}}\right) ,\forall x\in K,\forall
\varepsilon \leq \varepsilon _{0}\text{ \ ,}
\end{equation}
and $\mathcal{N}^{s}\left( \Omega \right) $ is the space of $\left(
f_{\varepsilon }\right) _{\varepsilon }\in C^{\infty }\left( \Omega \right)
^{\left] 0,1\right[ }$ satisfying for every compact $K$ of $\Omega $,$%
\forall \alpha \in \mathbb{Z}_{+}^{n},\exists C>0,\forall k>0,$ $\exists
\varepsilon _{0}\in \left] 0,1\right[ ,$
\begin{equation}
\left| \partial ^{\alpha }f_{\varepsilon }\left( x\right) \right| \leq C\exp
\left( -k\varepsilon ^{-\frac{1}{2s-1}}\right) ,\forall x\in K,\forall
\varepsilon \leq \varepsilon _{0}
\end{equation}
We show that $\mathcal{G}^{s}\left( \Omega \right) $\ contains the space of
Gevrey ultradistibutions of order $(3s-1),$ and the following diagram of
embeddings is commutative

\begin{equation}
\begin{array}{ccc}
D^{s}(\Omega ) & \rightarrow & \mathcal{G}^{s}(\Omega ) \\
& \searrow & \uparrow \\
&  & E_{3s-1}^{\prime }(\Omega )
\end{array}
\end{equation}

The Gevrey microlocal analysis in the framework of the algebra $\mathcal{G}%
^{s}(\Omega )$ consists first in introducing the algebra of regular
generalized Gevrey ultradistributions $\mathcal{G}^{s,\infty }(\Omega )$\
and the proof of the following fundamental result
\begin{equation}
\mathcal{G}^{s,\infty }(\Omega )\cap E_{3s-1}^{\prime }(\Omega
)=D^{s}(\Omega )
\end{equation}
Then, we define the generalized Gevrey wave front of $f\in \mathcal{G}%
^{s}\left( \Omega \right) $, denoted $WF_{g}^{s}\left( f\right) $, and give
its main properties.

Finally, we give an application of this generalized Gevrey microlocal
analysis. The product of two generalized Gevrey ultradistributions always
exists, but there is no final description of the generalized wave front of
this product. This problem is also still posed in the Colombeau algebra of
generalized functions. In \cite{Hor-Kun}, the well-known H\"{o}rmander's
result on the wave front of the product of two distributions, has been
extended to the case of two Colombeau generalized functions. We show this
result in the case of two generalized Gevrey ultradistributions, namely we
obtain the following theorem.

\begin{theorem}
Let $f,g\in \mathcal{G}^{s}\left( \Omega \right) $, satisfying $\forall x\in
\Omega ,$%
\begin{equation*}
\left( x,0\right) \notin WF_{g}^{s}\left( f\right) +WF_{g}^{s}\left(
g\right) ,
\end{equation*}
then
\begin{equation*}
WF_{g}^{s}\left( fg\right) \subseteq \left( WF_{g}^{s}\left( f\right)
+WF_{g}^{s}\left( g\right) \right) \cup WF_{g}^{s}\left( f\right) \cup
WF_{g}^{s}\left( g\right)
\end{equation*}
\end{theorem}

\section{Generalized Gevrey ultrdistributions}

According to the construction of Colombeau algebras of generalized
functions, we introduce an algebra of moderate elements and its ideal of
null elements depending on the Gevrey order $s>1.$

\begin{definition}
The space of moderate elements, denoted $\mathcal{E}_{m}^{s}\left( \Omega
\right) ,$ is the space of $\left( f_{\varepsilon }\right) _{\varepsilon
}\in C^{\infty }\left( \Omega \right) ^{\left] 0,1\right[ }$ satisfying for
every compact subset $K$ of $\Omega $, $\forall \alpha \in \mathbb{Z}%
_{+}^{n},\exists C>0,\exists k>0,$ $\exists \varepsilon _{0}\in \left] 0,1%
\right[ ,$ such that
\begin{equation}
\left| \partial ^{\alpha }f_{\varepsilon }\left( x\right) \right| \leq C\exp
\left( k\varepsilon ^{-\frac{1}{2s-1}}\right) ,\forall x\in K,\forall
\varepsilon \leq \varepsilon _{0}  \label{1*1}
\end{equation}
The space of null elements, denoted $\mathcal{N}^{s}\left( \Omega \right) ,$
is the space of $\left( f_{\varepsilon }\right) _{\varepsilon }\in C^{\infty
}\left( \Omega \right) ^{\left] 0,1\right[ }$ satisfying for every compact $%
K $ of $\Omega $,$\forall \alpha \in \mathbb{Z}_{+}^{n},\exists C>0,\forall
k>0,$ $\exists \varepsilon _{0}\in \left] 0,1\right[ ,$ such that
\begin{equation}
\left| \partial ^{\alpha }f_{\varepsilon }\left( x\right) \right| \leq C\exp
\left( -k\varepsilon ^{-\frac{1}{2s-1}}\right) ,\forall x\in K,\forall
\varepsilon \leq \varepsilon _{0}  \label{1*2}
\end{equation}
\end{definition}

The main properties of the spaces $\mathcal{E}_{m}^{s}\left( \Omega \right) $
and $\mathcal{N}^{s}\left( \Omega \right) $\ are given in the following
proposition.

\begin{proposition}
1) The space of moderate elements $\mathcal{E}_{m}^{s}\left( \Omega \right) $
is an algebra stable by derivation.

2) The space $\mathcal{N}^{s}\left( \Omega \right) $ is an ideal of $%
\mathcal{E}_{m}^{s}\left( \Omega \right) .$
\end{proposition}

\begin{proof}
1) Let $\left( f_{\varepsilon }\right) _{\varepsilon },\left( g_{\varepsilon
}\right) _{\varepsilon }\in \mathcal{E}_{m}^{s}\left( \Omega \right) $ and $%
K $ be a compact of $\Omega $, then $\forall \beta \in \mathbb{Z}_{+}^{n},$ $%
\exists C_{1}=C_{1}\left( \beta \right) >0,\exists k_{1}=k_{1}\left( \beta
\right) >0,\exists \varepsilon _{1\beta }\in \left] 0,1\right[ $ such that $%
\forall x\in K,\forall \varepsilon \leq \varepsilon _{1\beta },$%
\begin{equation}
\left| \partial ^{\beta }f_{\varepsilon }\left( x\right) \right| \leq
C_{1}\exp \left( k_{1}\varepsilon ^{-\frac{1}{2s-1}}\right)  \label{4}
\end{equation}
$\forall \beta \in \mathbb{Z}_{+}^{n},$ $\exists C_{2}=C_{2}\left( \beta
\right) >0,\exists k_{2}=K_{2}\left( \beta \right) >0,\exists \varepsilon
_{2\beta }\in \left] 0,1\right[ $, such that $\forall x\in K,\forall
\varepsilon \leq \varepsilon _{2\beta },$%
\begin{equation}
\left| \partial ^{\beta }g_{\varepsilon }\left( x\right) \right| \leq
C_{2}\exp \left( k_{2}\varepsilon ^{-\frac{1}{2s-1}}\right)  \label{5}
\end{equation}
Let $\alpha \in \mathbb{Z}_{+}^{n},$ then
\begin{equation*}
\left| \partial ^{\alpha }\left( f_{\varepsilon }g_{\varepsilon }\right)
\left( x\right) \right| \leq \sum_{\beta =0}^{\alpha }\binom{\alpha }{\beta }%
\left| \partial ^{\alpha -\beta }f_{\varepsilon }\left( x\right) \right|
\left| \partial ^{\beta }g_{\varepsilon }\left( x\right) \right|
\end{equation*}
For $k=\max \left\{ k_{1}\left( \beta \right) :\beta \leq \alpha \right\}
+\max \left\{ k_{2}\left( \beta \right) :\beta \leq \alpha \right\}
,\varepsilon \leq \min \left\{ \varepsilon _{1\beta },\varepsilon _{2\beta
};\left| \beta \right| \leq \left| \alpha \right| \right\} $ and $x\in K,$
we have
\begin{eqnarray*}
\exp \left( -k\varepsilon ^{-\frac{1}{2s-1}}\right) \left| \partial ^{\alpha
}\left( f_{\varepsilon }g_{\varepsilon }\right) \left( x\right) \right|
&\leq &\sum_{\beta =0}^{\alpha }\binom{\alpha }{\beta }\exp \left(
-k_{1}\varepsilon ^{-\frac{1}{2s-1}}\right) \left| \partial ^{\alpha -\beta
}f_{\varepsilon }\left( x\right) \right| \times \\
&&\times \exp \left( -k_{2}\varepsilon ^{-\frac{1}{2s-1}}\right) \left|
\partial ^{\beta }g_{\varepsilon }\left( x\right) \right| \\
&\leq &\sum_{\beta =0}^{\alpha }\binom{\alpha }{\beta }C_{1}\left( \alpha
-\beta \right) C_{2}\left( \beta \right) =C\left( \alpha \right) ,
\end{eqnarray*}
i.e. $\left( f_{\varepsilon }g_{\varepsilon }\right) _{\varepsilon }\in
\mathcal{E}_{m}^{s}\left( \Omega \right) $.

It is clear, from (\ref{4}) that for every compact $K$ of $\Omega $, $%
\forall \beta \in \mathbb{Z}_{+}^{n},$ $\exists C_{1}=C_{1}\left( \beta
+1\right) >0,\exists k_{1}=k_{1}\left( \beta +1\right) >0,\exists
\varepsilon _{1\beta }\in \left] 0,1\right[ $ such that $\forall x\in
K,\forall \varepsilon \leq \varepsilon _{1\beta },$%
\begin{equation*}
\left| \partial ^{\beta }\left( \partial f_{\varepsilon }\right) \left(
x\right) \right| \leq C_{1}\exp \left( k_{1}\varepsilon ^{-\frac{1}{2s-1}%
}\right) ,
\end{equation*}
i.e. $\left( \partial f_{\varepsilon }\right) _{\varepsilon }\in \mathcal{E}%
_{m}^{s}\left( \Omega \right) .$

2) If $\left( g_{\varepsilon }\right) _{\varepsilon }\in \mathcal{N}%
^{s}\left( \Omega \right) ,$ $\ $for every $K$ compact of $\Omega ,$ $%
\forall \beta \in \mathbb{Z}_{+}^{m},\exists C_{2}=C_{2}\left( \beta \right)
>0,\forall k_{2}>0,\exists \varepsilon _{2\beta }\in \left] 0,1\right[ ,$%
\begin{equation}
\left| \partial ^{\alpha }g_{\varepsilon }\left( x\right) \right| \leq
C_{2}\exp \left( -k_{2}\varepsilon ^{-\frac{1}{2s-1}}\right) ,\forall x\in
K,\forall \varepsilon \leq \varepsilon _{2\beta }
\end{equation}

Let $\alpha \in \mathbb{Z}_{+}^{m}$ and $k>0,$ then
\begin{eqnarray*}
\exp \left( k\varepsilon ^{-\frac{1}{2s-1}}\right) \left| \partial ^{\alpha
}\left( f_{\varepsilon }g_{\varepsilon }\right) \left( x\right) \right|
&\leq &\exp \left( k\varepsilon ^{-\frac{1}{2s-1}}\right) \sum_{\beta
=0}^{\alpha }\binom{\alpha }{\beta }\left| \partial ^{\alpha -\beta
}f_{\varepsilon }\left( x\right) \right| \times \\
&&\times \left| \partial ^{\beta }g_{\varepsilon }\left( x\right) \right|
\end{eqnarray*}
Let $k_{2}=\max \left\{ k_{1}\left( \alpha -\beta \right) :\beta \leq \alpha
\right\} +k$ and $\varepsilon \leq \min \left\{ \varepsilon _{1\beta
},\varepsilon _{2\beta };\beta \leq \alpha \right\} ,$ then $\forall x\in K$%
,
\begin{eqnarray*}
\exp \left( k\varepsilon ^{-\frac{1}{2s-1}}\right) \left| \partial ^{\alpha
}\left( f_{\varepsilon }g_{\varepsilon }\right) \left( x\right) \right|
&\leq &\sum_{\beta =0}^{\alpha }\binom{\alpha }{\beta }\left[ \exp \left(
-k_{1}\varepsilon ^{-\frac{1}{2s-1}}\right) \left| \partial ^{\alpha -\beta
}f_{\varepsilon }\left( x\right) \right| \right. \\
&&\left. \exp \left( k_{2}\varepsilon ^{-\frac{1}{2s-1}}\right) \left|
\partial ^{\beta }g_{\varepsilon }\left( x\right) \right| \right] \\
&\leq &\sum_{\beta =0}^{\alpha }\binom{\alpha }{\beta }C_{1}\left( \alpha
-\beta \right) C_{2}\left( \beta \right) =C\left( \alpha \right) ,
\end{eqnarray*}
which shows that $\left( f_{\varepsilon }g_{\varepsilon }\right)
_{\varepsilon }\in \mathcal{N}^{s}\left( \Omega \right) $
\end{proof}

\begin{remark}
The algebra of moderate elements $\mathcal{E}_{m}^{s}\left( \Omega \right) $%
\ is not necessary stable by $s-$ultradifferentiable operators.
\end{remark}

\begin{definition}
The algebra of generalized Gevrey ultradistributions of order $s>1$, denoted
$\mathcal{G}^{s}\left( \Omega \right) ,$ is the quotient algebra
\begin{equation}
\mathcal{G}^{s}\left( \Omega \right) =\frac{\mathcal{E}_{m}^{s}\left( \Omega
\right) }{\mathcal{N}^{s}\left( \Omega \right) }
\end{equation}
\end{definition}

\begin{remark}
We have $\mathcal{N}^{s}\left( \Omega \right) \subset \mathcal{N}\left(
\Omega \right) \subset \mathcal{E}_{m}\left( \Omega \right) \subset \mathcal{%
E}_{m}^{s}\left( \Omega \right) ,$ where $\mathcal{N}\left( \Omega \right) \
$is the Colombeau algebra of null elements and $\mathcal{E}_{m}\left( \Omega
\right) $ the Colombeau algebra of moderate elements.
\end{remark}

Our construction permits also to define the Gevrey type algebra of
generalized complex numbers, denoted $\mathcal{C}^{s},$ by
\begin{equation*}
\mathcal{C}^{s}=\frac{\mathcal{E}_{0}^{s}}{\mathcal{N}_{0}^{s}}\text{ \ ,}
\end{equation*}
where
\begin{equation*}
\mathcal{E}_{0}^{s}=\left\{ \left( a_{\varepsilon }\right) _{\varepsilon
}\in \mathbb{C}^{\left] 0,1\right[ };\exists k>0\text{, }\left|
a_{\varepsilon }\right| =O\left( \exp \left( k\varepsilon ^{-\frac{1}{2s-1}%
}\right) \right) \right\}
\end{equation*}
and
\begin{equation*}
\mathcal{N}_{0}^{s}=\left\{ \left( a_{\varepsilon }\right) _{\varepsilon
}\in \mathbb{C}^{\left] 0,1\right[ };\forall k>0\text{, }\left|
a_{\varepsilon }\right| =O\left( \exp \left( -k\varepsilon ^{-\frac{1}{2s-1}%
}\right) \right) \right\}
\end{equation*}
Some local properties of the algebra $\mathcal{G}^{s}\left( \Omega \right) $%
\ can be studied thanks to the algebra of generalized complex numbers $%
\mathcal{C}^{s}.$

\section{Embedding of Gevrey ultradistributions}

We recall some definitions and results on Gevrey ultradistributions, see
\cite{Kom}, \cite{LM} or \cite{Rod}.

\begin{definition}
A function $f$ $\in E^{s}\left( \Omega \right) $ if $f\in C^{\infty }\left(
\Omega \right) $ and for every compact subset $K$ of $\Omega $, $\exists
c>0,\forall \alpha \in \mathbb{Z}_{+}^{m},$%
\begin{equation}
\sup_{x\in K}\left| \partial ^{\alpha }f\left( x\right) \right| \leq
c^{\left| \alpha \right| +1}\left( \alpha !\right) ^{s}
\end{equation}
\end{definition}

Obviously we have $E^{t}\left( \Omega \right) \subset E^{s}\left( \Omega
\right) $ if $1\leq t\leq s.$ It is well known that $E^{1}\left( \Omega
\right) =A\left( \Omega \right) $ is the space of all real analytic
functions in $\Omega $ and if we denote by $D^{s}\left( \Omega \right) $ the
space $E^{s}\left( \Omega \right) \cap C_{0}^{\infty }\left( \Omega \right)
, $ then $D^{s}\left( \Omega \right) $\ is non trivial if and only if $s>1.$
The topological dual of $D^{s}\left( \Omega \right) ,$ denoted $%
D_{s}^{\prime }\left( \Omega \right) ,$\ is called the space of Gevrey
ultradistributions of order $s.$\ The space $E_{s}^{\prime }\left( \Omega
\right) $ is the topological dual of $E^{s}\left( \Omega \right) $\ and is
identified with\ the space of Gevrey ultradistributions with compact
supports.

\begin{definition}
A differential operator of infinite order $P\left( D\right)
=\sum\limits_{\gamma \in \mathbb{Z}_{+}^{m}}a_{\gamma }D^{\gamma }$ is
called an ultradifferential operator of class $s$ or $s$-ultradifferential
operator, if for every $h>0$ there exists $c>0$ such that $\forall \gamma
\in \mathbb{Z}_{+}^{m},$
\begin{equation}
\left| a_{\gamma }\right| \leq c\frac{h^{\left| \gamma \right| }}{\left(
\gamma !\right) ^{s}}  \label{12}
\end{equation}
\end{definition}

The importance of $s$-ultradifferential operator\ is in the following result.

\begin{proposition}
Let $T\in E_{s}^{\prime }\left( \Omega \right) $ and ${supp}T\subset K,$
then there exists an $s$-ultradifferential operator $P\left( D\right)
=\sum\limits_{\gamma \in \mathbb{Z}_{+}^{m}}a_{\gamma }D^{\gamma }$ , $M>0$\
and continuous functions $\ f_{\gamma }\in C_{0}\left( K\right) $ such that $%
\sup\limits_{\gamma \in \mathbb{Z}_{+}^{m},x\in \mathbb{R}^{m}}\left|
f_{\gamma }\left( x\right) \right| \leq M$ and
\begin{equation}
T=\sum_{\gamma \in \mathbb{Z}_{+}^{m}}a_{\gamma }D^{\gamma }f_{\gamma }
\end{equation}
\end{proposition}

The space $\mathcal{S}^{\left( s\right) },s>1,$ see \cite{Gel}, is the space
of functions $\varphi \in C^{\infty }\left( \Omega \right) $ such that $%
\forall b>0,$ we have
\begin{equation}
\sigma _{b,s}\left( \varphi \right) =\sup_{\alpha ,\beta \in \mathbb{Z}%
_{+}^{m}}\int \frac{\left| x\right| ^{\left| \beta \right| }}{b^{\left|
\alpha +\beta \right| }\alpha !^{s}\beta !^{s}}\left| \partial ^{\alpha
}\varphi \left( x\right) \right| dx<\infty  \label{2-1}
\end{equation}

\begin{lemma}
There exists a$\ \phi \in \mathcal{S}^{\left( s\right) }$ satisfying
\begin{equation*}
\int \phi \left( x\right) dx=1\text{ and }\int x^{\alpha }\phi \left(
x\right) dx=0,\forall \alpha \in \mathbb{Z}_{+}^{m}\backslash \{0\}
\end{equation*}
\end{lemma}

\begin{proof}
For an example of function $\phi \in \mathcal{S}^{\left( s\right) }$
satisfying these conditions, take the\ Fourier transform of a function of
the class $D^{\left( s\right) }\left( \Omega \right) $ equals $1$ in
neighborhood of the origin. Here $D^{\left( s\right) }\left( \Omega \right) $
denotes the projective Gevrey space of order $s,$ i. e. $D^{\left( s\right)
}\left( \Omega \right) =E^{(s)}\left( \Omega \right) \cap C_{0}^{\infty
}\left( \Omega \right) ,$ where $f$ $\in E^{(s)}\left( \Omega \right) ,$\ if
$f\in C^{\infty }\left( \Omega \right) $ and for every compact subset $K$ of
$\Omega $,$\forall b>0,$ $\exists c>0,\forall \alpha \in \mathbb{Z}_{+}^{m},$%
\begin{equation}
\sup_{x\in K}\left| \partial ^{\alpha }f\left( x\right) \right| \leq
cb^{\left| \alpha \right| }\left( \alpha !\right) ^{s}
\end{equation}
\end{proof}

\begin{definition}
The net $\phi _{\varepsilon }=\varepsilon ^{-m}\phi \left( ./\varepsilon
\right) ,\varepsilon \in \left] 0,1\right[ ,$ where $\phi $\ satisfies the
Lemma, is called a net of mollifiers.
\end{definition}

The space $E^{s}\left( \Omega \right) $ is embedded into $\mathcal{G}%
^{s}\left( \Omega \right) $ by the standard canonical injection
\begin{equation*}
\begin{array}{ccc}
I: & E^{s}\left( \Omega \right) \rightarrow & \mathcal{G}^{s}\left( \Omega
\right) \\
& f\rightarrow & \left[ f\right] =cl\left( f_{\varepsilon }\right) ,
\end{array}
\end{equation*}
where $f_{\varepsilon }=f$ , $\forall \varepsilon \in \left] 0,1\right[ $.

The following proposition gives the natural embedding of Gevrey
ultradistributions into $\mathcal{G}^{s}\left( \Omega \right) .$

\begin{theorem}
\label{pro-inj}The map
\begin{equation}
\begin{array}{ccc}
J: & E_{3s-1}^{\prime }\left( \Omega \right) \rightarrow & \mathcal{G}%
^{s}\left( \Omega \right) \\
& T\rightarrow & \left[ T\right] =cl\left( \left( T\ast \phi _{\varepsilon
}\right) _{/\Omega }\right) _{\varepsilon }
\end{array}
\end{equation}
is an embedding$.$
\end{theorem}

\begin{proof}
Let $T\in E_{3s-1}^{\prime }\left( \Omega \right) $ with ${supp}T\subset K,$
then there exists an $\left( 3s-1\right) $-ultradifferential operator $%
P\left( D\right) =\sum\limits_{\gamma \in \mathbb{Z}_{+}^{m}}a_{\gamma
}D^{\gamma }$ and continuous functions $f_{\gamma }$ with ${supp}f_{\gamma
}\subset K,\forall \gamma \in \mathbb{Z}_{+}^{m},$ and $\sup\limits_{\gamma
\in \mathbb{Z}_{+}^{m},x\in K}\left| f_{\gamma }\left( x\right) \right| \leq
M,$ such that
\begin{equation}
T=\sum_{\gamma \in \mathbb{Z}_{+}^{m}}a_{\gamma }D^{\gamma }f_{\gamma }
\end{equation}
We have
\begin{equation*}
T\ast \phi _{\varepsilon }\left( x\right) =\sum_{\gamma \in \mathbb{Z}%
_{+}^{m}}a_{\gamma }\frac{\left( -1\right) ^{\left| \gamma \right| }}{%
\varepsilon ^{\left| \gamma \right| }}\int f_{\gamma }\left( x+\varepsilon
y\right) D^{\gamma }\phi \left( y\right) dy
\end{equation*}
Let $\alpha \in \mathbb{Z}_{+}^{m},$ then
\begin{equation*}
\left| \partial ^{\alpha }\left( T\ast \phi _{\varepsilon }\left( x\right)
\right) \right| \leq \sum_{\gamma \in \mathbb{Z}_{+}^{m}}a_{\gamma }\frac{1}{%
\varepsilon ^{\left| \gamma +\alpha \right| }}\int \left| f_{\gamma }\left(
x+\varepsilon y\right) \right| \left| D^{\gamma +\alpha }\phi \left(
y\right) \right| dy
\end{equation*}
From (\ref{12}) and the inequality
\begin{equation}
\left( \beta +\alpha \right) !^{t}\leq 2^{t\left| \beta +\alpha \right|
}\alpha !^{t}\beta !^{t},\text{ }\forall t\geq 1,  \label{6}
\end{equation}
we have, $\forall h>0,\exists c>0,$ such that
\begin{eqnarray*}
\left| \partial ^{\alpha }\left( T\ast \phi _{\varepsilon }\left( x\right)
\right) \right| &\leq &\sum_{\gamma \in \mathbb{Z}_{+}^{m}}c\frac{h^{\left|
\gamma \right| }}{\gamma !^{3s-1}}\frac{1}{\varepsilon ^{\left| \gamma
+\alpha \right| }}\int \left| f_{\gamma }\left( x+\varepsilon y\right)
\right| \left| D^{\gamma +\alpha }\phi \left( y\right) \right| dy \\
&\leq &\sum_{\gamma \in \mathbb{Z}_{+}^{m}}c\alpha !^{3s-1}\frac{2^{\left(
3s-1\right) \left| \gamma +\alpha \right| }h^{\left| \gamma \right| }}{%
\left( \gamma +\alpha \right) !^{2s-1}}\frac{1}{\varepsilon ^{\left| \gamma
+\alpha \right| }}b^{\left| \gamma +\alpha \right| }\times \\
&&\qquad \times \int \left| f_{\gamma }\left( x+\varepsilon y\right) \right|
\frac{\left| D^{\gamma +\alpha }\phi \left( y\right) \right| }{b^{\left|
\gamma +\alpha \right| }\left( \gamma +\alpha \right) !^{s}}dy,
\end{eqnarray*}
then for $h>\frac{1}{2},$
\begin{eqnarray*}
\frac{1}{\alpha !^{3s-1}}\left| \partial ^{\alpha }\left( T\ast \phi
_{\varepsilon }\left( x\right) \right) \right| &\leq &\sigma _{b,s}\left(
\phi \right) Mc\sum_{\gamma \in \mathbb{Z}_{+}^{m}}2^{-\left| \gamma \right|
}\frac{\left( 2^{3s}bh\right) ^{\left| \gamma +\alpha \right| }}{\left(
\gamma +\alpha \right) !^{2s-1}}\frac{1}{\varepsilon ^{\left| \gamma +\alpha
\right| }} \\
&\leq &C\exp \left( k_{1}\varepsilon ^{-\frac{1}{2s-1}}\right) ,
\end{eqnarray*}
i.e.
\begin{equation}
\left| \partial ^{\alpha }\left( T\ast \phi _{\varepsilon }\left( x\right)
\right) \right| \leq C\left( \alpha \right) \exp \left( k_{1}\varepsilon ^{-%
\frac{1}{2s-1}}\right) ,  \label{1-2}
\end{equation}
where $k_{1}=\left( 2s-1\right) \left( 2^{3s}bh\right) ^{\frac{1}{2s-1}}$.

Suppose that $\left( T\ast \phi _{\varepsilon }\right) _{\varepsilon }\in
\mathcal{N}^{s}\left( \Omega \right) ,$ then for every compact $L$ of $%
\Omega ,$ $\exists C>0,$ $\forall k>0,$ $\exists \varepsilon _{0}\in \left]
0,1\right[ ,$%
\begin{equation}
\left| T\ast \phi _{\varepsilon }\left( x\right) \right| \leq C\exp \left(
-k\varepsilon ^{-\frac{1}{2s-1}}\right) ,\forall x\in L,\forall \varepsilon
\leq \varepsilon _{0}  \label{7}
\end{equation}
Let $\chi \in D^{3s-1}\left( \Omega \right) $ and $\chi =1$ in neighborhood
of $K$, then $\forall \psi \in E^{3s-1}\left( \Omega \right) ,$%
\begin{equation*}
\left\langle T,\psi \right\rangle =\left\langle T,\chi \psi \right\rangle
=\lim_{\varepsilon \rightarrow 0}\int \left( T\ast \phi _{\varepsilon
}\right) \left( x\right) \chi \left( x\right) \psi \left( x\right) dx
\end{equation*}
Consequently, from (\ref{7}), we obtain
\begin{equation*}
\left| \int \left( T\ast \phi _{\varepsilon }\right) \left( x\right) \chi
\left( x\right) \psi \left( x\right) dx\right| \leq C\exp \left(
-k\varepsilon ^{-\frac{1}{2s-1}}\right) ,\forall \varepsilon \leq
\varepsilon _{0},
\end{equation*}
which gives $\left\langle T,\psi \right\rangle =0$
\end{proof}

\begin{remark}
We have $C\left( \alpha \right) =\alpha !^{3s-1}\sigma _{b,s}\left( \phi
\right) Mc$ in (\ref{1-2}).
\end{remark}

In order to show that the following diagram of embeddings
\begin{equation*}
\begin{array}{ccc}
D^{s}\left( \Omega \right) & \rightarrow & \mathcal{G}^{s}\left( \Omega
\right) \\
& \searrow & \uparrow \\
&  & E_{3s-1}^{\prime }\left( \Omega \right)
\end{array}
\end{equation*}
is commutative, we have to prove the following fundamental result.

\begin{proposition}
\label{pro2}Let $f\in D^{s}\left( \Omega \right) $ and $\left( \phi
_{\varepsilon }\right) _{\varepsilon }$ be a net of mollifiers, then
\begin{equation*}
\left( f-\left( f\ast \phi _{\varepsilon }\right) _{/\Omega }\right)
_{\varepsilon }\in \mathcal{N}^{s}\left( \Omega \right)
\end{equation*}
\end{proposition}

\begin{proof}
Let $f\in D^{s}\left( \Omega \right) ,$ then there exists a constant $C>0,$
such that
\begin{equation}
\left| \partial ^{\alpha }f\left( x\right) \right| \leq C^{\left| \alpha
\right| +1}\alpha !^{s},\forall \alpha \in \mathbb{Z}_{+}^{m},\forall x\in
\Omega
\end{equation}
Let $\alpha \in \mathbb{Z}_{+}^{m}$, the Taylor formula and the properties
of $\phi _{\varepsilon }$ give
\begin{equation*}
\partial ^{\alpha }\left( f\ast \phi _{\varepsilon }-f\right) \left(
x\right) =\sum_{\left| \beta \right| =N}\int \frac{\left( \varepsilon
y\right) ^{\beta }}{\beta !}\partial ^{\alpha +\beta }f\left( \xi \right)
\phi \left( y\right) dy,
\end{equation*}
where $x\leq \xi \leq x+\varepsilon y.$ Consequently, for $b>0$, we have
\begin{eqnarray*}
\left| \partial ^{\alpha }\left( f\ast \phi _{\varepsilon }-f\right) \left(
x\right) \right| &\leq &\varepsilon ^{N}\sum_{\left| \beta \right| =N}\int
\frac{\left| y\right| ^{N}}{\beta !}\left| \partial ^{\alpha +\beta }f\left(
\xi \right) \right| \left| \phi \left( y\right) \right| dy \\
&\leq &\alpha !^{s}\varepsilon ^{N}\sum_{\left| \beta \right| =N}\beta
!^{2s-1}2^{s\left| \alpha +\beta \right| }b^{\left| \beta \right| }\int
\frac{\left| \partial ^{\alpha +\beta }f\left( \xi \right) \right| }{\left(
\alpha +\beta \right) !^{s}}\times \\
&&\times \frac{\left| y\right| ^{\left| \beta \right| }}{b^{\left| \beta
\right| }\beta !^{s}}\left| \phi \left( y\right) \right| dy
\end{eqnarray*}
Let $k>0$ and $T>0,$ then
\begin{eqnarray*}
\left| \partial ^{\alpha }\left( f\ast \phi _{\varepsilon }-f\right) \left(
x\right) \right| &\leq &\alpha !^{s}\left( \varepsilon N^{2s-1}\right)
^{N}\left( k^{2s-1}T\right) ^{-N}\times \\
&&\times \sum_{\left| \beta \right| =N}\int 2^{s\left| \alpha +\beta \right|
}\left( k^{2s-1}bT\right) ^{\left| \beta \right| }\frac{\left| \partial
^{\alpha +\beta }f\left( \xi \right) \right| }{\left( \alpha +\beta \right)
!^{s}}\times \\
&&\times \frac{\left| y\right| ^{\left| \beta \right| }}{b^{\left| \beta
\right| }\beta !^{s}}\left| \phi \left( y\right) \right| dy \\
&\leq &\alpha !^{s}\left( \varepsilon N^{2s-1}\right) ^{N}\left(
k^{2s-1}T\right) ^{-N}\times \\
&&\times C\sigma _{b,s}\left( \phi \right) \left( 2^{s}C\right) ^{\left|
\alpha \right| }\sum_{\left| \beta \right| =N}\left( 2^{s}k^{2s-1}bT\right)
^{\left| \beta \right| }C^{\left| \beta \right| },
\end{eqnarray*}
hence, taking $2^{s}k^{2s-1}bTC\leq \frac{1}{2a},$ with $a>1$, we obtain
\begin{eqnarray}
\left| \partial ^{\alpha }\left( f\ast \phi _{\varepsilon }-f\right) \left(
x\right) \right| &\leq &\alpha !^{s}\left( \varepsilon N^{2s-1}\right)
^{N}\left( k^{2s-1}T\right) ^{-N}\times  \notag \\
&&\times C\sigma _{b,s}\left( \phi \right) \left( 2^{s}C\right) ^{\left|
\alpha \right| }a^{-N}\sum_{\left| \beta \right| =N}\left( \frac{1}{2}%
\right) ^{\left| \beta \right| }  \label{1*5} \\
&\leq &\sigma _{b,s}\left( \phi \right) C^{\left| \alpha \right| +1}\alpha
!^{s}\left( \varepsilon N^{2s-1}\right) ^{N}\left( k^{2s-1}T\right)
^{-N}a^{-N}  \notag
\end{eqnarray}
Let $\varepsilon _{0}\in \left] 0,1\right[ $ such that $\varepsilon _{0}^{%
\frac{1}{2s-1}}\frac{\ln a}{k}<1$ and take $T>2^{2s-1},$ then
\begin{equation*}
\left( T^{\frac{1}{2s-1}}-1\right) >1>\frac{\ln a}{k}\varepsilon ^{\frac{1}{%
2s-1}},\forall \varepsilon \leq \varepsilon _{0},
\end{equation*}
in particular, we have
\begin{equation*}
\left( \frac{\ln a}{k}\varepsilon ^{\frac{1}{2s-1}}\right) ^{-1}T^{\frac{1}{%
2s-1}}-\left( \frac{\ln a}{k}\varepsilon ^{\frac{1}{2s-1}}\right) ^{-1}>1
\end{equation*}
Then, there exists $N=N\left( \varepsilon \right) \in \mathbb{Z}^{+},$ such
that
\begin{equation*}
\left( \frac{\ln a}{k}\varepsilon ^{\frac{1}{2s-1}}\right) ^{-1}<N<\left(
\frac{\ln a}{k}\varepsilon ^{\frac{1}{2s-1}}\right) ^{-1}T^{\frac{1}{2s-1}},
\end{equation*}
i.e.
\begin{equation}
1\leq \frac{\ln a}{k}\varepsilon ^{\frac{1}{2s-1}}N\leq T^{\frac{1}{2s-1}},
\label{2-2}
\end{equation}
which gives
\begin{equation*}
a^{-N}\leq \exp \left( -k\varepsilon ^{-\frac{1}{2s-1}}\right) \text{ \ \ \
\ \ and\ \ \ \ \ }\frac{\varepsilon N^{2s-1}}{k^{2s-1}T}\leq \left( \frac{1}{%
\ln a}\right) ^{2s-1}<1,
\end{equation*}
if we choose $\ln a>1$. Finally, from (\ref{1*5}), we have
\begin{equation}
\left| \partial ^{\alpha }\left( f\ast \phi _{\varepsilon }-f\right) \left(
x\right) \right| \leq C\exp \left( -k\varepsilon ^{-\frac{1}{2s-1}}\right) ,
\label{1*4}
\end{equation}
i.e. $f\ast \phi _{\varepsilon }-f\in \mathcal{N}^{s}\left( \Omega \right) $
\end{proof}

From the proof, see (\ref{1*5}), we obtained in fact the following result. \

\begin{corollary}
Let $f\in D^{s}\left( \Omega \right) $, then for every compact $K$ of $%
\Omega ,\exists C>0,\forall \alpha \in \mathbb{Z}_{+}^{m},\forall
k>0,\exists \varepsilon _{0}\in \left] 0,1\right[ ,\forall x\in K,\forall
\varepsilon \leq \varepsilon _{0}$ ,
\begin{equation}
\left| \partial ^{\alpha }\left( f\ast \phi _{\varepsilon }-f\right) \left(
x\right) \right| \leq C^{\left| \alpha \right| +1}\alpha !^{s}\exp \left(
-k\varepsilon ^{-\frac{1}{2s-1}}\right)  \label{1*3}
\end{equation}
\end{corollary}

Let $\Omega ^{\prime }$ be an open subset of $\Omega $ and let $f=\left(
f_{\varepsilon }\right) _{\varepsilon }+\mathcal{N}^{s}\left( \Omega \right)
\in \mathcal{G}^{s}\left( \Omega \right) $, the restriction of $f$ to $%
\Omega ^{\prime },$ denoted $f_{/\Omega ^{\prime }}$, is defined as
\begin{equation}
\left( f_{\varepsilon /\Omega ^{\prime }}\right) _{\varepsilon }+\mathcal{N}%
^{s}\left( \Omega ^{\prime }\right) \in \mathcal{G}^{s}\left( \Omega
^{\prime }\right)
\end{equation}
One can easily prove that the functor $\Omega \rightarrow \mathcal{G}%
^{s}\left( \Omega \right) $ defines a presheaf in the same way as in the
case of the algebra of Colombeau generalized functions $\mathcal{G}\left(
\Omega \right) $, see for example \cite{GKOS}. Consequently, we introduce
the support of $f\in $ $\mathcal{G}^{s}\left( \Omega \right) $, denoted ${%
supp}_{g}f,$ as the complement of the largest open set $U$ such that $%
f_{/U}=0$.

\begin{definition}
The space of elements of $\mathcal{G}^{s}\left( \Omega \right) $\ with
compact supports is denoted $\mathcal{G}_{c}^{s}\left( \Omega \right) .$
\end{definition}

It is not difficult to prove the following result.

\begin{proposition}
\label{ref3}1) The space $\mathcal{G}_{c}^{s}\left( \Omega \right) $ is the
space of $f=cl\left( f_{\varepsilon }\right) _{\varepsilon }\in \mathcal{G}%
^{s}\left( \Omega \right) $ satisfying, $\exists K$ a compact subset of $%
\Omega ,\exists \varepsilon _{0}\in \left] 0,1\right[ $,$\forall \varepsilon
\in \left] 0,\varepsilon _{0}\right[ ,{supp}f_{\varepsilon }\subset K.$
\newline
2) Let $f=cl\left( f_{\varepsilon }\right) _{\varepsilon }\in \mathcal{G}%
_{c}^{s}\left( \Omega \right) ,$ then $\exists C>0,\exists k_{1}>0,\exists
\varepsilon _{0}>0,\forall k_{2}>0,\forall \xi \in \mathbb{R}^{m},\forall
\varepsilon \leq \varepsilon _{0},$
\begin{equation}
{\left| \widehat{f_{\varepsilon }}\left( \xi \right) \right| \leq C\exp
\left( k_{1}\varepsilon ^{-\frac{1}{2s-1}}+k_{2}\left| \xi \right| ^{\frac{1%
}{s}}\right) }
\end{equation}
\end{proposition}

\section{Regular generalized Gevrey ultradistributions}

To develop a local or microlocal analysis with respect to a ''good space of
regular elements'' one needs first to define these regular elements, then
the notion of singular support and its microlocalisation.

\begin{definition}
We define $\mathcal{E}_{m}^{s,\infty }\left( \Omega \right) $ the space of
regular elements as the space of $\left( f_{\varepsilon }\right)
_{\varepsilon }\in \left( C^{\infty }\left( \Omega \right) \right) ^{\left]
0,1\right[ }$ satisfying, for every compact subset $K$ of $\Omega $, $%
\exists C>0,\exists k>0,$ $\exists \varepsilon _{0}\in \left] 0,1\right[
,\forall \alpha \in \mathbb{Z}_{+}^{n},\forall x\in K,\forall \varepsilon
\leq \varepsilon _{0}$,
\begin{equation}
\left| \partial ^{\alpha }f_{\varepsilon }\left( x\right) \right| \leq
C^{\left| \alpha \right| +1}\alpha !^{s}\exp \left( k\varepsilon ^{-\frac{1}{%
2s-1}}\right)
\end{equation}
\end{definition}

\begin{proposition}
1) The space $\mathcal{E}_{m}^{s,\infty }\left( \Omega \right) $ is an
algebra stable by the action of $s$-ultradifferential operators.

2) The space $\mathcal{N}_{\ast }^{s}\left( \Omega \right) :=\mathcal{N}%
^{s}\left( \Omega \right) \cap \mathcal{E}_{m}^{s,\infty }\left( \Omega
\right) $ is an ideal of $\mathcal{E}_{m}^{s,\infty }\left( \Omega \right) .$
\end{proposition}

\begin{proof}
1) Let $\left( f_{\varepsilon }\right) _{\varepsilon },\left( g_{\varepsilon
}\right) _{\varepsilon }\in \mathcal{E}_{m}^{s,\infty }\left( \Omega \right)
$ and $K$ be a compact of $\Omega $, then $\exists c_{1}>0,\exists
k_{1}>0,\exists \varepsilon _{1}\in \left] 0,1\right[ $ such that $\forall
x\in K,\forall \alpha \in \mathbb{Z}_{+}^{m},\forall \varepsilon \leq
\varepsilon _{1\beta },$%
\begin{equation}
\left| \partial ^{\alpha }f_{\varepsilon }\left( x\right) \right| \leq
c_{1}^{\left| \alpha \right| +1}{\alpha !^{s}}\exp \left( k_{1}\varepsilon
^{-\frac{1}{2s-1}}\right)
\end{equation}

We have also $\exists c_{2}>0,\exists k_{2}>0,\exists \varepsilon _{2}\in %
\left] 0,1\right[ $ such that $\forall x\in K,\forall \alpha \in \mathbb{Z}%
_{+}^{m},\forall \varepsilon \leq \varepsilon _{2},$%
\begin{equation}
\left| \partial ^{\alpha }g_{\varepsilon }\left( x\right) \right| \leq
c_{2}^{\left| \alpha \right| +1}{\alpha !^{s}}\exp \left( k_{2}\varepsilon
^{-\frac{1}{2s-1}}\right)
\end{equation}

Let $\alpha \in \mathbb{Z}_{+}^{n}$, then
\begin{equation*}
\frac{1}{\alpha !^{s}}\left| \partial ^{\alpha }\left( f_{\varepsilon
}g_{\varepsilon }\right) \left( x\right) \right| \leq \sum_{\beta
=0}^{\alpha }\binom{\alpha }{\beta }\frac{1}{\left( \alpha -\beta \right)
!^{s}}\left| \partial ^{\alpha -\beta }f_{\varepsilon }\left( x\right)
\right| \frac{1}{\beta !^{s}}\left| \partial ^{\beta }g_{\varepsilon }\left(
x\right) \right|
\end{equation*}
Let $\varepsilon \leq \min \left( \varepsilon _{1},\varepsilon _{2}\right) $
and $k=k_{1}+k_{2},$ then we have $\forall \alpha \in \mathbb{Z}%
_{+}^{m},\forall x\in K,$%
\begin{eqnarray*}
\exp \left( -k\varepsilon ^{-\frac{1}{2s-1}}\right) \frac{1}{\alpha !^{s}}%
\left| \partial ^{\alpha }\left( f_{\varepsilon }g_{\varepsilon }\right)
\left( x\right) \right| &\leq &\sum_{\beta =0}^{\alpha }\binom{\alpha }{%
\beta }\frac{\exp \left( -k_{1}\varepsilon ^{-\frac{1}{2s-1}}\right) }{%
\left( \alpha -\beta \right) !^{s}}\left| \partial ^{\alpha -\beta
}f_{\varepsilon }\left( x\right) \right| \times \\
&&\times \frac{\exp \left( -k_{2}\varepsilon ^{-\frac{1}{2s-1}}\right) }{%
\beta !^{s}}\left| \partial ^{\beta }g_{\varepsilon }\left( x\right) \right|
\\
&\leq &\sum_{\beta =0}^{\alpha }\binom{\alpha }{\beta }c_{1}^{\left| \alpha
-\beta \right| }c_{2}^{\left| \beta \right| } \\
&\leq &2^{^{\left| \alpha \right| }}\left( c_{1}+c_{2}\right) ^{\left|
\alpha \right| },
\end{eqnarray*}
i. e. $\left( f_{\varepsilon }\right) _{\varepsilon }\left( g_{\varepsilon
}\right) _{\varepsilon }\in \mathcal{E}_{m}^{s,\infty }\left( \Omega \right)
.$

Let now $P\left( D\right) =\sum a_{\gamma }D^{\gamma }$ be an$\ s$%
-ultradifferential operator, then $\forall h>0,\exists b>0,$ such that
\begin{eqnarray*}
\exp \left( k_{1}\varepsilon ^{-\frac{1}{2s-1}}\right) \frac{1}{\alpha !^{s}}%
\left| \partial ^{\alpha }\left( P\left( D\right) f_{\varepsilon }\left(
x\right) \right) \right| &\leq &\exp \left( k_{1}\varepsilon ^{-\frac{1}{2s-1%
}}\right) \sum_{\gamma \in \mathbb{Z}_{+}^{m}}b\frac{h^{\left| \gamma
\right| }}{\gamma !^{s}}\frac{1}{\alpha !^{s}}\left| \partial ^{\alpha
+\gamma }f_{\varepsilon }\left( x\right) \right| \\
&\leq &b\exp \left( k_{1}\varepsilon ^{-\frac{1}{2s-1}}\right) \sum_{\gamma
\in \mathbb{Z}_{+}^{m}}\frac{2^{s\left| \alpha +\gamma \right| }h^{\left|
\gamma \right| }}{\left( \alpha +\gamma \right) !^{s}}\left| \partial
^{\alpha +\gamma }f_{\varepsilon }\left( x\right) \right| \\
&\leq &b\sum_{\gamma \in \mathbb{Z}_{+}^{m}}2^{s\left| \alpha +\gamma
\right| }h^{\left| \gamma \right| }c_{1}^{\left| \alpha +\gamma \right| },
\end{eqnarray*}
hence, for $2^{s}hc_{1}\leq \frac{1}{2},$ we have
\begin{equation*}
\exp \left( k\varepsilon ^{-\frac{1}{2s-1}}\right) \frac{1}{\alpha !^{s}}%
\left| \partial ^{\alpha }\left( P\left( D\right) f_{\varepsilon }\left(
x\right) \right) \right| \leq c^{\prime }\left( 2^{s}c_{1}\right) ^{\left|
\alpha \right| },
\end{equation*}
which shows that $\left( P\left( D\right) f_{\varepsilon }\right)
_{\varepsilon }\in \mathcal{E}_{m}^{s,\infty }\left( \Omega \right) .$

2) The fact that $\mathcal{N}_{\ast }^{s}\left( \Omega \right) =\mathcal{N}%
^{s}\left( \Omega \right) \cap \mathcal{E}_{m}^{s,\infty }\left( \Omega
\right) \subset \mathcal{E}_{m}^{s}\left( \Omega \right) $ and $\mathcal{N}%
^{s}\left( \Omega \right) $ is an ideal of $\mathcal{E}_{m}^{s}\left( \Omega
\right) ,$ then $\mathcal{N}_{\ast }^{s}\left( \Omega \right) $ is an ideal
of $\mathcal{E}_{m}^{s,\infty }\left( \Omega \right) $
\end{proof}

\begin{remark}
If the inclusion $\mathcal{N}^{s}\left( \Omega \right) \subset \mathcal{E}%
_{m}^{s,\infty }\left( \Omega \right) $ holds, then $\mathcal{N}_{\ast
}^{s}\left( \Omega \right) =\mathcal{N}^{s}\left( \Omega \right) .$
\end{remark}

Now, we define the Gevrey regular elements of $\mathcal{G}^{s}\left( \Omega
\right) .$

\begin{definition}
The algebra of regular generalized Gevrey ultradistributions of order $s>1$,
denoted $\mathcal{G}^{s,\infty }\left( \Omega \right) ,$ is the quotient
algebra
\begin{equation*}
\mathcal{G}^{s,\infty }\left( \Omega \right) =\frac{\mathcal{E}%
_{m}^{s,\infty }\left( \Omega \right) }{\mathcal{N}_{\ast }^{s}\left( \Omega
\right) }
\end{equation*}
\end{definition}

\begin{remark}
It is clear that $E^{s}\left( \Omega \right) \hookrightarrow \mathcal{G}%
^{s,\infty }\left( \Omega \right) .$
\end{remark}

\begin{definition}
We define the $\mathcal{G}^{s,\infty }$-singular support of a generalized
Gevrey ultradistribution $f$ $\in \mathcal{G}^{s}\left( \Omega \right) ,$
denoted $s$-$singsupp_{g}\left( f\right) ,$ as the complement of the largest
open set $\Omega ^{\prime }$ such that $f\in \mathcal{G}^{s,\infty }\left(
\Omega ^{\prime }\right) .$
\end{definition}

\begin{proposition}
\label{pro1}Let $f=cl\left( f_{\varepsilon }\right) _{\varepsilon }\in
\mathcal{G}_{c}^{s}\left( \Omega \right) ,$ then $f$ is regular if and only
if $\exists k_{1}>0,\exists k_{2}>0,\exists C>0,\exists \varepsilon
_{1}>0,\forall \varepsilon \leq \varepsilon _{1}$ , such that
\begin{equation}
\left| \widehat{f_{\varepsilon }}\left( \xi \right) \right| \leq C\exp
\left( k_{1}\varepsilon ^{-\frac{1}{2s-1}}-k_{2}\left| \xi \right| ^{\frac{1%
}{s}}\right) ,\forall \xi \in \mathbb{R}^{m}  \label{3-2}
\end{equation}
\end{proposition}

\begin{proof}
Suppose that $f=cl\left( f_{\varepsilon }\right) _{\varepsilon }\in \mathcal{%
G}_{c}^{s}\left( \Omega \right) \cap \mathcal{G}^{s,\infty }\left( \Omega
\right) ,$ then $\exists C_{1}>0,\exists k_{1}>0,\exists \varepsilon _{1}>0,$
$\forall \alpha \in \mathbb{Z}_{+}^{n},\forall x\in K,\forall \varepsilon
\leq \varepsilon _{1},$ ${supp}f_{\varepsilon }\subset K,$ such that
\begin{equation}
\left| \partial ^{\alpha }f_{\varepsilon }\right| \leq C_{1}^{\left| \alpha
\right| +1}\alpha !^{s}\exp \left( k_{1}\varepsilon ^{-\frac{1}{2s-1}}\right)
\end{equation}
Consequently we have, $\forall \alpha \in \mathbb{Z}_{+}^{m},$
\begin{equation*}
\left| \xi ^{\alpha }\right| \left| \widehat{f_{\varepsilon }}\left( \xi
\right) \right| \leq \left| \int \exp \left( -ix\xi \right) \partial
^{\alpha }f_{\varepsilon }\left( x\right) dx\right| ,
\end{equation*}
then, $\exists C>0,\forall \varepsilon \leq \varepsilon _{1},$%
\begin{equation*}
\left| \xi \right| ^{\left| \alpha \right| }\left| \widehat{f_{\varepsilon }}%
\left( \xi \right) \right| \leq C^{\left| \alpha \right| +1}\alpha !^{s}\exp
\left( k_{1}\varepsilon ^{-\frac{1}{2s-1}}\right)
\end{equation*}
For $\alpha \in \mathbb{Z}_{+}^{m},\exists N\in \mathbb{Z}_{+}$ such that
\begin{equation*}
\frac{N}{s}\leq \left| \alpha \right| <\frac{N}{s}+1,
\end{equation*}
so
\begin{eqnarray*}
\left| \xi \right| ^{\frac{N}{s}}\left| \widehat{f_{\varepsilon }}\left( \xi
\right) \right| &\leq &C^{\left| \alpha \right| +1}\left| \alpha \right|
^{\left| \alpha \right| s}\exp \left( k_{1}\varepsilon ^{-\frac{1}{2s-1}%
}\right) \\
&\leq &C^{N+1}N^{N}\exp \left( k_{1}\varepsilon ^{-\frac{1}{2s-1}}\right)
\end{eqnarray*}
Hence $\exists C>0,\forall N\in \mathbb{Z}^{+},$%
\begin{equation*}
\left| \widehat{f_{\varepsilon }}\left( \xi \right) \right| \leq
C^{N+1}\left| \xi \right| ^{-\frac{N}{s}}N!\exp \left( k_{1}\varepsilon ^{-%
\frac{1}{2s-1}}\right) ,
\end{equation*}
which gives
\begin{equation*}
\left| \widehat{f_{\varepsilon }}\left( \xi \right) \right| \exp \left(
\frac{1}{2C}\left| \xi \right| ^{\frac{1}{s}}\right) \leq C\exp \left(
k_{1}\varepsilon ^{-\frac{1}{2s-1}}\right) \sum 2^{-N},
\end{equation*}
or
\begin{equation*}
\left| \widehat{f_{\varepsilon }}\left( \xi \right) \right| \leq C^{\prime
}\exp \left( k_{1}\varepsilon ^{-\frac{1}{2s-1}}-\frac{1}{2C}\left| \xi
\right| ^{\frac{1}{s}}\right) ,
\end{equation*}
i.e. we have (\ref{3-2}).

Suppose now that (\ref{3-2}) is valid, then $\forall \varepsilon \leq
\varepsilon _{0},$
\begin{equation*}
\left| \partial ^{\alpha }f_{\varepsilon }\left( x\right) \right| \leq
C_{1}\exp \left( k_{1}\varepsilon ^{-\frac{1}{2s-1}}\right) \int \left| \xi
^{\alpha }\right| \exp \left( -k_{2}\left| \xi \right| ^{\frac{1}{s}}\right)
d\xi
\end{equation*}
Due to the inequality $t^{N}\leq N!\exp \left( t\right) ,\forall t>0,$ then $%
\exists C_{2}=C\left( k_{2}\right) $ such that
\begin{equation*}
\left| \xi ^{\alpha }\right| \exp \left( -\frac{k_{2}}{2}\left| \xi \right|
^{\frac{1}{s}}\right) \leq C_{2}^{\left| \alpha \right| }\alpha !^{s},
\end{equation*}
then
\begin{eqnarray*}
\left| \partial ^{\alpha }f_{\varepsilon }\left( x\right) \right| &\leq
&C_{1}\exp \left( k_{1}\varepsilon ^{-\frac{1}{2s-1}}\right) C_{2}^{\left|
\alpha \right| }\alpha !^{s}\int \exp \left( -\frac{k_{2}}{2}\left| \xi
\right| ^{\frac{1}{s}}\right) d\xi \\
&\leq &C^{\left| \alpha \right| +1}\alpha !^{s}\exp \left( k_{1}\varepsilon
^{-\frac{1}{2s-1}}\right) ,\text{ }
\end{eqnarray*}
where $C=\max \left( C_{1}\int \exp \left( -\frac{k_{2}}{2}\left| \xi
\right| ^{\frac{1}{s}}\right) d\xi ,C_{2}\right) $, i.e. $f\in \mathcal{G}%
^{s,\infty }\left( \Omega \right) $
\end{proof}

The algebra $\mathcal{G}^{s,\infty }\left( \Omega \right) $ plays the same
role as the Oberguggenberger subalgebra of regular elements $\mathcal{G}%
^{\infty }\left( \Omega \right) $ of the Colombeau algebra $\mathcal{G}%
\left( \Omega \right) $, see \cite{Ober}$.$

\begin{theorem}
We have
\begin{equation*}
\mathcal{G}^{s,\infty }\left( \Omega \right) \cap E_{3s-1}^{\prime }\left(
\Omega \right) =D^{s}\left( \Omega \right)
\end{equation*}
\end{theorem}

\begin{proof}
Let $T\in \mathcal{G}^{s,\infty }\left( \Omega \right) \cap E_{3s-1}^{\prime
}\left( \Omega \right) ,$ with ${supp}T=K$ and $\phi _{\varepsilon }$ be a
net of mollifiers with $\check{\phi}=\phi $ and let $\chi \in D^{s}\left(
\Omega \right) $ such that $\chi =1$ on $K.$ As $\left[ T\right] \in
\mathcal{G}^{s,\infty }\left( \Omega \right) ,$ then $\exists
c_{1}>0,\exists k_{1}>0,\exists k_{2}>0,\exists \varepsilon _{1}>0,\forall
\varepsilon \leq \varepsilon _{1},$%
\begin{equation*}
\left| \widehat{\chi \left( T\ast \phi _{\varepsilon }\right) }\left( \xi
\right) \right| \leq c_{1}e^{k_{1}\varepsilon ^{-\frac{1}{2s-1}}-k_{2}\left|
\xi \right| ^{\frac{1}{s}}},
\end{equation*}
then
\begin{eqnarray*}
\left| \widehat{\chi \left( T\ast \phi _{\varepsilon }\right) }\left( \xi
\right) -\widehat{T}\left( \xi \right) \right| &=&\left| \widehat{\chi
\left( T\ast \phi _{\varepsilon }\right) }\left( \xi \right) -\widehat{\chi T%
}\left( \xi \right) \right| \\
&=&\left| \left\langle T,\left( \chi e^{-i\xi .}\right) \ast \phi
_{\varepsilon }-\left( \chi e^{-i\xi .}\right) \right\rangle \right|
\end{eqnarray*}
As $E_{3s-1}^{\prime }\left( \Omega \right) \subset E_{s}^{\prime }\left(
\Omega \right) ,$ then $\exists L$ a compact subset of $\Omega $ such that $%
\forall h>0,\exists c>0,$ and
\begin{equation*}
\left| \widehat{\chi \left( T\ast \phi _{\varepsilon }\right) }\left( \xi
\right) -\widehat{T}\left( \xi \right) \right| \leq c\sup_{\alpha \in
\mathbb{Z}_{+}^{m},x\in L}\frac{h^{\left| \alpha \right| }}{\alpha !^{s}}%
\left| \left( \partial ^{\alpha }\left( \chi e^{-i\xi }\ast \phi
_{\varepsilon }-\chi e^{-i\xi .}\right) \left( x\right) \right) \right|
\end{equation*}
We have $e^{-i\xi }\chi \in D^{s}\left( \Omega \right) ,$ from the corollary
3-9, $\exists c_{2}>0,\forall k_{3}>0,\exists \eta >0,\forall \varepsilon
\leq \eta ,$%
\begin{equation*}
\sup_{\alpha \in \mathbb{Z}_{+}^{m},x\in L}\frac{c_{2}^{\left| \alpha
\right| }}{\alpha !^{s}}\left| \partial ^{\alpha }\left( \chi e^{-i\xi
.}\ast \phi _{\varepsilon }-\chi e^{-i\xi .}\right) \left( x\right) \right|
\leq c_{2}e^{-k_{3}\varepsilon ^{-\frac{1}{2s-1}}},
\end{equation*}
so there exists $c^{\prime }>0$, such that
\begin{equation*}
\left| \widehat{T}\left( \xi \right) -\widehat{\chi \left( T\ast \phi
_{\varepsilon }\right) }\left( \xi \right) \right| \leq c^{\prime
}e^{-k_{3}\varepsilon ^{-\frac{1}{2s-1}}}
\end{equation*}
Let $\varepsilon \leq \min \left( \eta ,\varepsilon _{1}\right) ,$ then
\begin{eqnarray*}
\left| \widehat{T}\left( \xi \right) \right| &\leq &\left| \widehat{T}\left(
\xi \right) -\widehat{\chi \left( T\ast \phi _{\varepsilon }\right) }\left(
\xi \right) \right| +\left| \widehat{\chi \left( T\ast \phi _{\varepsilon
}\right) }\left( \xi \right) \right| \\
&\leq &c^{\prime }e^{-k_{3}\varepsilon ^{-\frac{1}{2s-1}}}+c_{1}e^{k_{1}%
\varepsilon ^{-\frac{1}{2s-1}}-k_{2}\left| \xi \right| ^{\frac{1}{s}}}
\end{eqnarray*}
Take $c=\max \left( c^{\prime },c_{1}\right) ,$ $\varepsilon =\left( \dfrac{%
k_{1}}{\left( k_{2}-r\right) \left| \xi \right| ^{\frac{1}{s}}}\right)
^{2s-1},r\in \left] 0,k_{2}\right[ $ and $k_{3}=\dfrac{k_{1}r}{k_{2}-r}$,
then $\exists \delta =2r>0,\exists c>0$ such that
\begin{equation*}
\left| \widehat{T}\left( \xi \right) \right| \leq ce^{-\delta \left| \xi
\right| ^{\frac{1}{s}}},
\end{equation*}
which means $T\in E^{s}\left( \Omega \right) $. As ${supp}T=K$ is a compact
then we have $T\in D^{s}\left( \Omega \right) $
\end{proof}

\section{Generalized Gevrey wave front}

The defined local regularity and its Fourier characterization, studied in
the previous section, allow us to define the fundamental concept of every
microlocal analysis, i. e. the generalized Gevrey wave front of a
generalized Gevrey ultradistribution.

\begin{definition}
We define $\sum_{g}^{s}\left( f\right) \subset \mathbb{R}^{m}\backslash
\left\{ 0\right\} ,f\in \mathcal{G}_{c}^{s}\left( \Omega \right) $, as the
complement of the set of points having a conic neighborhood $\Gamma $ such
that $\exists k_{1}>0,\exists k_{2}>0,\exists C>0,\exists \varepsilon
_{0}\in \left] 0,1\right[ ,\forall \xi \in \Gamma ,\forall \varepsilon \leq
\varepsilon _{0},$
\begin{equation}
\left| \widehat{f_{\varepsilon }}\left( \xi \right) \right| \leq c\exp
\left( k_{1}\varepsilon ^{-\frac{1}{2s-1}}-k_{2}\left| \xi \right| ^{\frac{1%
}{s}}\right)  \label{3-1}
\end{equation}
\end{definition}

\begin{proposition}
\label{ref5}For every $f\in \mathcal{G}_{c}^{s}\left( \Omega \right) $, we
have

1) The set $\sum_{g}^{s}\left( f\right) $ is a close cone.

2) $\sum_{g}^{s}\left( f\right) =\emptyset \Longleftrightarrow f\in \mathcal{%
G}^{s,\infty }\left( \Omega \right) .$

3) $\sum_{g}^{s}\left( \psi f\right) \subset \sum_{g}^{s}\left( f\right)
,\forall \psi \in E^{s}\left( \Omega \right) .$
\end{proposition}

\begin{proof}
The proofs of 1) is easy, 2) holds from the proposition \ref{pro1}. Let us
proof 3), if $\xi _{0}\notin \sum_{g}^{s}\left( f\right) $, then $\exists
\Gamma $ a conic neighborhood of $\xi _{0},\exists k_{1}>0,\exists
k_{2}>0,\exists c_{1}>0,\exists \varepsilon _{0}\in \left] 0,1\right[ ,$ $%
\forall \xi \in \Gamma ,\forall \varepsilon \leq \varepsilon _{0},$
\begin{equation}
\left| \widehat{f_{\varepsilon }}\left( \xi \right) \right| \leq c_{1}\exp
\left( k_{1}\varepsilon ^{-\frac{1}{2s-1}}-k_{2}\left| \xi \right| ^{\frac{1%
}{s}}\right)  \label{3-02}
\end{equation}
Let $\chi \in D^{s}\left( \Omega \right) ,$ $\chi =1$ on neighborhood of $%
suppf$, so $\chi \psi \in D^{s}\left( \Omega \right) $, hence $\exists
k_{3}>0,\exists c_{2}>0,\forall \xi \in \mathbb{R}^{m},$
\begin{equation}
\left| \widehat{\chi \psi }\left( \xi \right) \right| \leq c_{2}\exp \left(
-k_{3}\left| \xi \right| ^{\frac{1}{s}}\right) ,  \label{3-3}
\end{equation}
Let $\Lambda $ be a conic neighborhood of $\xi _{0}$ such that, $\overline{%
\Lambda }\subset \Gamma ,$ we have, for a fixed $\xi \in \Lambda $,
\begin{equation*}
\widehat{\chi \psi f_{\varepsilon }}\left( \xi \right) =\int_{A}\widehat{%
f_{\varepsilon }}\left( \eta \right) \widehat{\chi \psi }\left( \eta -\xi
\right) d\eta +\int_{B}\widehat{f_{\varepsilon }}\left( \eta \right)
\widehat{\chi \psi }\left( \eta -\xi \right) d\eta ,
\end{equation*}
where $A=\left\{ \eta ;\left| \xi -\eta \right| ^{\frac{1}{s}}\leq \delta
\left( \left| \xi \right| ^{\frac{1}{s}}+\left| \eta \right| ^{\frac{1}{s}%
}\right) \right\} $; $B=\left\{ \eta ;\left| \xi -\eta \right| ^{\frac{1}{s}%
}>\delta \left( \left| \xi \right| ^{\frac{1}{s}}+\left| \eta \right| ^{%
\frac{1}{s}}\right) \right\} $. We choose $\delta $ sufficiently small such
that $A\subset \Gamma $ and $\dfrac{\left| \xi \right| }{2^{s}}<\left| \eta
\right| <2^{s}\left| \xi \right| $. Then $\exists \varepsilon _{0}\in \left]
0,1\right[ ,\forall \varepsilon \leq \varepsilon _{0},$
\begin{eqnarray}
\left| \int_{A}\widehat{f_{\varepsilon }}\left( \eta \right) \widehat{\chi
\psi }\left( \eta -\xi \right) d\eta \right| &\leq &c_{1}c_{2}\exp \left(
k_{1}\varepsilon ^{-\frac{1}{2s-1}}-\frac{k_{2}}{2}\left| \xi \right| ^{%
\frac{1}{s}}\right) \times  \notag \\
&&\times \left| \int_{A}\exp \left( -k_{3}\left| \eta -\xi \right| ^{\frac{1%
}{s}}\right) d\eta \right|  \notag \\
&\leq &c^{\prime }\exp \left( k_{1}\varepsilon ^{-\frac{1}{2s-1}}-\frac{k_{2}%
}{2}\left| \xi \right| ^{\frac{1}{s}}\right)  \label{3-4}
\end{eqnarray}
As $f\in \mathcal{G}_{c}^{s}\left( \Omega \right) ,$ from \ proposition \ref
{ref3}, $\exists c>0,\exists \mu _{1}>0,\exists \varepsilon _{1}\in \left]
0,1\right[ ,\forall \mu _{2}>0,\forall \xi \in \mathbb{R}^{m},\forall
\varepsilon \leq \varepsilon _{1},$ such that
\begin{equation*}
\left| \widehat{f_{\varepsilon }}\left( \xi \right) \right| \leq c\exp
\left( \mu _{1}\varepsilon ^{-\frac{1}{2s-1}}+\mu _{2}\left| \xi \right| ^{%
\frac{1}{s}}\right) ,
\end{equation*}
hence, for $\varepsilon \leq \min \left( \varepsilon _{0},\varepsilon
_{1}\right) ,$ we have
\begin{eqnarray*}
\left| \int_{B}\widehat{f_{\varepsilon }}\left( \eta \right) \widehat{\psi }%
\left( \eta -\xi \right) d\eta \right| &\leq &cc_{2}\exp \left( \mu
_{1}\varepsilon ^{-\frac{1}{2s-1}}\right) \int_{B}\exp \left( \mu _{2}\left|
\eta \right| ^{\frac{1}{s}}-k_{3}\left| \eta -\xi \right| ^{\frac{1}{s}%
}\right) d\eta \\
&\leq &c^{\prime \prime }\exp \left( \mu _{1}\varepsilon ^{-\frac{1}{2s-1}%
}\right) \int_{B}\exp \left( \mu _{2}\left| \eta \right| ^{\frac{1}{s}%
}-k_{3}\delta \left( \left| \xi \right| ^{\frac{1}{s}}+\left| \eta \right| ^{%
\frac{1}{s}}\right) \right) d\eta ,
\end{eqnarray*}
then, taking $\mu _{2}<k_{3}\delta ,$ we obtain
\begin{equation}
\left| \int_{B}\widehat{f_{\varepsilon }}\left( \eta \right) \widehat{\psi }%
\left( \eta -\xi \right) d\eta \right| \leq c\exp \left( \mu _{1}\varepsilon
^{-\frac{1}{2s-1}}-k_{3}\delta \left| \xi \right| ^{\frac{1}{s}}\right)
\label{3-5}
\end{equation}
Consequently, (\ref{3-4}) and (\ref{3-5}) give $\xi _{0}\notin
\sum_{g}^{s}\left( \psi f\right) $
\end{proof}

\begin{definition}
Let $f\in \mathcal{G}^{s}\left( \Omega \right) $ and $x_{0}\in \Omega $, the
cone of $s$-singular directions of $f$ at $x_{0}$, denoted $%
\sum_{g,x_{0}}^{s}\left( f\right) $, is
\begin{equation}
\sum\nolimits_{g,x_{0}}^{s}\left( f\right) =\bigcap \left\{
\sum\nolimits_{g}^{s}\left( \phi f\right) :\phi \in D^{s}\left( \Omega
\right) \text{ and }\phi =1\text{ on a neighborhood of }x_{0}\right\}
\label{3-6}
\end{equation}
\end{definition}

\begin{lemma}
\label{lem1}Let $f\in \mathcal{G}^{s}\left( \Omega \right) $, then
\begin{equation*}
\sum\nolimits_{g,x_{0}}^{s}\left( f\right) =\emptyset \Longleftrightarrow
x_{0}\notin s\text{-}singsupp_{g}\left( f\right)
\end{equation*}
\end{lemma}

\begin{proof}
Let $x_{0}\notin s$-$singsupp_{g}\left( f\right) ,$ i.e. $\exists U\subset
\Omega $ an open neighborhood of $x_{0}$ such that $f\in \mathcal{G}%
^{s,\infty }\left( U\right) $, let $\phi \in D^{s}\left( U\right) $ such
that $\phi =1$ on a neighborhood of $x_{0},$ then $\phi f\in \mathcal{G}%
^{s,\infty }\left( \Omega \right) .$ Hence, from the proposition \ref{ref5},
$\sum\nolimits_{g}^{s}\left( \phi f\right) =\emptyset ,$ \ i.e. $%
\sum\nolimits_{g,x_{0}}^{s}\left( f\right) =\emptyset .$

Suppose now $\sum\nolimits_{g,x_{0}}^{s}\left( f\right) =\emptyset ,$ let $%
r>0$ such that $B\left( 0,2r\right) \subset \Omega $ and let $\psi \in
D^{s}\left( B\left( 0,2r\right) \right) $ such that $0\leq $ $\psi \leq 1$
and $\psi =1$ on $B\left( 0,r\right) $. Let $\psi _{j}\left( x\right) =\psi
\left( 3^{j}\left( x-x_{0}\right) \right) $ then it is clear that $%
supp\left( \psi _{j}\right) \subset B\left( x_{0},\frac{2r}{3^{j}}\right)
\subset \Omega $ and $\psi _{j}=1$ on $B\left( x_{0},\frac{r}{3^{j}}\right)
, $ we have $\forall \phi \in D^{s}\left( \Omega \right) $ with $\phi =1$ on
a neighborhood $U$ of $x_{0},$ $\exists j\in \mathbb{Z}^{+}$ such that $%
supp\left( \psi _{j}\right) \subset U$, then $\psi _{j}f_{\varepsilon }=\psi
_{j}\phi f_{\varepsilon }$ and from proposition \ref{ref5}, we have
\begin{equation*}
\sum\nolimits_{g}^{s}\left( \psi _{j}f\right) \subset
\sum\nolimits_{g}^{s}\left( \phi f\right) ,
\end{equation*}
which gives
\begin{equation}
\bigcap\limits_{j\in \mathbb{Z}^{+}}\left( \sum\nolimits_{g}^{s}\left( \psi
_{j}f\right) \right) =\sum\nolimits_{g,x_{0}}^{s}\left( f\right) =\emptyset
\label{3-7}
\end{equation}
We have $\psi _{j}=1$ on $supp\left( \psi _{j+1}\right) ,$ then $%
\sum\nolimits_{g}^{s}\left( \psi _{j+1}f\right) \subset
\sum\nolimits_{g}^{s}\left( \psi _{j}f\right) $, so from (\ref{3-7}), there
exists $n\in \mathbb{Z}^{+}$ sufficiently large such that $\left( \psi
_{n}f\right) \in \mathcal{G}^{s,\infty }\left( \Omega \right) ,$ then $f\in
\mathcal{G}^{s,\infty }\left( B\left( x_{0},\frac{r}{3^{n}}\right) \right) ,$
which means. $x_{0}\notin s$-$singsupp_{g}\left( f\right) $
\end{proof}

Now, we are ready to give the definition of the generalized Gevrey wave
front and its main properties.

\begin{definition}
A point $\left( x_{0},\xi _{0}\right) \notin WF_{g}^{s}\left( f\right)
\subset \Omega \times \mathbb{R}^{m}\backslash \left\{ 0\right\} $ if $\xi
_{0}\notin \sum\nolimits_{g,x_{0}}^{s}\left( f\right) ,$ i.e. there exists $%
\phi \in D^{s}\left( \Omega \right) ,\phi \left( x\right) =1$ neighborhood
of $x_{0}$, and conic neighborhood $\Gamma $ of $\xi _{0}$, $\exists
k_{1}>0,\exists k_{2}>0,\exists c>0,\exists \varepsilon _{0}\in \left] 0,1%
\right[ ,$ such that $\forall \xi \in \Gamma ,\forall \varepsilon \leq
\varepsilon _{0},$
\begin{equation*}
\left| \widehat{\phi f_{\varepsilon }}\left( \xi \right) \right| \leq c\exp
\left( k_{1}\varepsilon ^{-\frac{1}{2s-1}}-k_{2}\left| \xi \right| ^{\frac{1%
}{s}}\right)
\end{equation*}
\end{definition}

The main properties of the generalized Gevrey wave front $WF_{g}^{s}$ are
given in the following proposition.

\begin{proposition}
Let $f\in \mathcal{G}^{s}\left( \Omega \right) $, then

1) The projection of $WF_{g}^{s}\left( f\right) $ on $\Omega $ is the $%
s-singsupp_{g}\left( f\right) $

2) If $f\in \mathcal{G}_{c}^{s}\left( \Omega \right) ,$ then the projection
of $WF_{g}^{s}\left( f\right) $ on $\mathbb{R}^{m}\backslash \left\{
0\right\} $ is $\sum_{g}^{s}\left( f\right) $

3) $\forall \alpha \in \mathbb{Z}_{+}^{m},WF_{g}^{s}\left( \partial ^{\alpha
}f\right) \subset WF_{g}^{s}\left( f\right) $

4) $\forall g\in \mathcal{G}^{s,\infty }\left( \Omega \right)
,WF_{g}^{s}\left( gf\right) \subset WF_{g}^{s}\left( f\right) $
\end{proposition}

\begin{proof}
1) and 2) holds from the definition, the proposition \ref{ref5} and lemma
\ref{lem1}. 3) Let $\left( x_{0},\xi _{0}\right) \notin WF_{g}^{s}\left(
f\right) $, then $\exists \phi \in D^{s}\left( \Omega \right) ,\phi \left(
x_{0}\right) =1$ on a neighborhood $\overline{U}$ of $x_{0}$, there exist a
conic neighborhood $\Gamma $ of $\xi
_{0},k_{1}>0,k_{2}>0,c_{1}>0,\varepsilon _{0}\in \left] 0,1\right[ ,$ such
that $\forall \xi \in \Gamma ,\forall \varepsilon \leq \varepsilon _{0},$
\begin{equation}
\left| \widehat{\phi f_{\varepsilon }}\left( \xi \right) \right| \leq
c_{1}\exp \left( k_{1}\varepsilon ^{-\frac{1}{2s-1}}-k_{2}\left| \xi \right|
^{\frac{1}{s}}\right)  \label{3-8}
\end{equation}
We have, for $\psi \in D^{s}\left( U\right) $ such that $\psi \left(
x_{0}\right) =1$,
\begin{eqnarray*}
\left| \widehat{\psi \partial f_{\varepsilon }}\left( \xi \right) \right|
&=&\left| \widehat{\partial \left( \psi f_{\varepsilon }\right) }\left( \xi
\right) -\widehat{\left( \partial \psi \right) f_{\varepsilon }}\left( \xi
\right) \right| \\
&\leq &\left| \xi \right| \left| \widehat{\psi \phi f_{\varepsilon }}\left(
\xi \right) \right| +\left| \widehat{\left( \partial \psi \right) \phi
f_{\varepsilon }}\left( \xi \right) \right|
\end{eqnarray*}
As $WF_{g}^{s}\left( \psi f\right) \subset WF_{g}^{s}\left( f\right) ,$ then
(\ref{3-8}) holds for both $\left| \widehat{\psi \phi f_{\varepsilon }}%
\left( \xi \right) \right| \ $and $\left| \widehat{\left( \partial \psi
\right) \phi f_{\varepsilon }}\left( \xi \right) \right| .$ So
\begin{eqnarray*}
\left| \xi \right| \left| \widehat{\psi \phi f_{\varepsilon }}\left( \xi
\right) \right| &\leq &c\left| \xi \right| \exp \left( k_{1}\varepsilon ^{-%
\frac{1}{2s-1}}-k_{2}\left| \xi \right| ^{\frac{1}{s}}\right) \\
&\leq &c^{\prime }\exp \left( k_{1}\varepsilon ^{-\frac{1}{2s-1}%
}-k_{3}\left| \xi \right| ^{\frac{1}{s}}\right) ,
\end{eqnarray*}
with $c^{\prime }>0,k_{3}>0$ such that $\left| \xi \right| \leq c^{\prime
}\exp \left( k_{2}-k_{3}\right) \left| \xi \right| ^{\frac{1}{s}}.$ Hence (%
\ref{3-8}) holds for $\left| \widehat{\psi \partial f_{\varepsilon }}\left(
\xi \right) \right| ,$ which proves $\left( x_{0},\xi _{0}\right) \notin
WF_{g}^{s}\left( \partial f\right) $

4) Let $\left( x_{0},\xi _{0}\right) \notin WF_{g}^{s}\left( f\right) $,
then $\exists \phi \in D^{s}\left( \Omega \right) ,\phi \left( x\right) =1$
on a neighborhood $U$ of $x_{0}$, there exist a conic neighborhood $\Gamma $
of $\xi _{0},$ $k_{1}>0,k_{2}>0,c_{1}>0,\varepsilon _{0}\in \left] 0,1\right[
,$ such that $\forall \xi \in \Gamma ,\forall \varepsilon \leq \varepsilon
_{0},$
\begin{equation*}
\left| \widehat{\phi f_{\varepsilon }}\left( \xi \right) \right| \leq
c_{1}\exp \left( k_{1}\varepsilon ^{-\frac{1}{2s-1}}-k_{2}\left| \xi \right|
^{\frac{1}{s}}\right)
\end{equation*}
Let $\psi \in D^{s}\left( \Omega \right) $ and $\psi =1$ on $supp\phi ,$
then $\widehat{\phi g_{\varepsilon }f_{\varepsilon }}=\widehat{\psi
g_{\varepsilon }}\ast \widehat{\phi f_{\varepsilon }}.$ We have $\psi g\in
\mathcal{G}^{s,\infty }\left( \Omega \right) ,$ then $\exists c_{2}>0,$ $%
\exists k_{3}>0,\exists k_{4}>0,\exists \varepsilon _{1}>0,\forall \xi \in
\mathbb{R}^{m},\forall \varepsilon \leq \varepsilon _{1},$
\begin{equation*}
\left| \widehat{\psi g_{\varepsilon }}\left( \xi \right) \right| \leq
c_{2}\exp \left( k_{3}\varepsilon ^{-\frac{1}{2s-1}}-k_{4}\left| \xi \right|
^{\frac{1}{s}}\right) ,
\end{equation*}
so
\begin{equation*}
\widehat{\phi g_{\varepsilon }f_{\varepsilon }}\left( \xi \right) =\int_{A}%
\widehat{\phi f_{\varepsilon }}\left( \eta \right) \widehat{\psi
g_{\varepsilon }}\left( \eta -\xi \right) d\eta +\int_{B}\widehat{\phi
f_{\varepsilon }}\left( \eta \right) \widehat{\psi g_{\varepsilon }}\left(
\eta -\xi \right) d\eta ,
\end{equation*}
where $A$ and $B$ are the same as in the proof of proposition \ref{ref5}. By
proposition \ref{ref3}, we have $\exists c>0,\exists \mu _{1}>0,\forall \mu
_{2}>0,\exists \varepsilon _{2}>0,\forall \xi \in \mathbb{R}^{m},\forall
\varepsilon \leq \varepsilon _{2},$%
\begin{equation*}
\left| \widehat{\phi f_{\varepsilon }}\left( \xi \right) \right| \leq c\exp
\left( \mu _{1}\varepsilon ^{-\frac{1}{2s-1}}+\mu _{2}\left| \xi \right| ^{%
\frac{1}{s}}\right)
\end{equation*}
The same steps as the proposition \ref{ref5}\ finish the proof
\end{proof}

\begin{corollary}
Let $P\left( x,D\right) =\sum\limits_{\left| \alpha \right| \leq m}a_{\alpha
}\left( x\right) D^{\alpha }$ be a partial differential operator with $%
\mathcal{G}^{s,\infty }\left( \Omega \right) $ coefficient, then
\begin{equation*}
WF_{g}^{s}\left( P\left( x,D\right) f\right) \subset WF_{g}^{s}\left(
f\right) ,\forall f\in \mathcal{G}^{s}\left( \Omega \right)
\end{equation*}
\end{corollary}

\begin{remark}
The reverse inclusion will give a generalized Gevrey microlocal
hypoellipticity of linear partial differential operators with generalized
coefficients. The case of generalized\ $\mathcal{G}^{\infty }-$microlocal
hypoellipticity has been studied recently in \cite{HorObPil}.
\end{remark}

\section{Generalized H\"{o}rmander's theorem}

We define $WF_{g}^{s}\left( f\right) +WF_{g}^{s}\left( g\right) ,$ where $%
f,g\in \mathcal{G}^{s}\left( \Omega \right) ,$ as the set
\begin{equation*}
\left\{ \left( x,\xi +\eta \right) ;\left( x,\xi \right) \in
WF_{g}^{s}\left( f\right) ,\left( x,\eta \right) \in WF_{g}^{s}\left(
g\right) \right\}
\end{equation*}

\begin{lemma}
\label{lem2}Let $\sum_{1}$, $\sum_{2}$ be closed cones in $\mathbb{R}%
^{m}\backslash \left\{ 0\right\} ,$ such that $0\notin \sum_{1}+\sum_{2}$ ,
then

i) $\overline{\sum\nolimits_{1}+\sum\nolimits_{2}}^{\mathbb{R}^{m}/\left\{
0\right\} }=\left( \sum\nolimits_{1}+\sum\nolimits_{2}\right) \cup
\sum\nolimits_{1}\cup \sum\nolimits_{2}$

ii) For any open conic neighborhood $\Gamma $ of $\sum\nolimits_{1}+\sum%
\nolimits_{2}$ in $\mathbb{R}^{m}\backslash \left\{ 0\right\} ,$ one can
find open conic neighborhoods of $\Gamma _{1},$ $\Gamma _{2}$ in $\mathbb{R}%
^{m}\backslash \left\{ 0\right\} $ of, respectively, $\sum\nolimits_{1},\sum%
\nolimits_{2}$ , such that
\begin{equation*}
\Gamma _{1}+\Gamma _{2}\subset \Gamma
\end{equation*}
\end{lemma}

\begin{proof}
See \cite{Hor-Kun}
\end{proof}

\begin{theorem}
Let $f,g\in \mathcal{G}^{s}\left( \Omega \right) $ , such that $\forall x\in
\Omega ,$%
\begin{equation}
\left( x,0\right) \notin WF_{g}^{s}\left( f\right) +WF_{g}^{s}\left(
g\right) ,  \label{4-1}
\end{equation}
then
\begin{equation}
WF_{g}^{s}\left( fg\right) \subseteq \left( WF_{g}^{s}\left( f\right)
+WF_{g}^{s}\left( g\right) \right) \cup WF_{g}^{s}\left( f\right) \cup
WF_{g}^{s}\left( g\right)  \label{4-2}
\end{equation}
\end{theorem}

\begin{proof}
Let $\left( x_{0},\xi _{0}\right) \notin \left( WF_{g}^{s}\left( f\right)
+WF_{g}^{s}\left( g\right) \right) \cup WF_{g}^{s}\left( f\right) \cup
WF_{g}^{s}\left( g\right) ,$ then $\exists \phi \in D^{s}\left( \Omega
\right) $, $\phi \left( x_{0}\right) =1,$ $\xi _{0}\notin \left(
\sum\nolimits_{g}^{s}\left( \phi f\right) +\sum\nolimits_{g}^{s}\left( \phi
g\right) \right) \cup \sum\nolimits_{g}^{s}\left( \phi f\right) \cup
\sum\nolimits_{g}^{s}\left( \phi g\right) =\overline{\sum\nolimits_{g}^{s}%
\left( \phi f\right) +\sum\nolimits_{g}^{s}\left( \phi g\right) }^{\mathbb{R}%
^{\ast }/0}.$ Let $\overline{\Gamma }_{0}$ be an open conic neighborhood of $%
\sum\nolimits_{g}^{s}\left( \phi f\right) +\sum\nolimits_{g}^{s}\left( \phi
g\right) $ in $\mathbb{R}^{m}\backslash \left\{ 0\right\} $ such that $\xi
_{0}\notin \overline{\Gamma }_{0}$. By lemma \ref{lem2} there exist open
cones $\Gamma _{1}$ and $\Gamma _{2}$ in $\mathbb{R}^{m}\backslash \left\{
0\right\} $ such that
\begin{equation*}
\sum\nolimits_{g}^{s}\left( \phi f\right) \subset \Gamma _{1},\text{ }%
\sum\nolimits_{g}^{s}\left( \phi g\right) \subset \Gamma _{2}\text{ and }%
\Gamma _{1}+\Gamma _{2}\subset \Gamma _{0}
\end{equation*}
Define $\Gamma =\mathbb{R}^{m}\backslash \overline{\Gamma }_{0},$ so
\begin{equation}
\Gamma \cap \Gamma _{2}=\emptyset \text{ and }\left( \Gamma -\Gamma
_{2}\right) \cap \Gamma _{1}=\emptyset  \label{4.3}
\end{equation}
Let $\xi \in \Gamma $ and $\varepsilon \in \left] 0,1\right[ $
\begin{eqnarray*}
\widehat{\phi f_{\varepsilon }\phi g_{\varepsilon }}\left( \xi \right)
&=&\left( \widehat{\phi f_{\varepsilon }}\ast \widehat{\phi g_{\varepsilon }}%
\right) \left( \xi \right) \\
&=&\int_{\Gamma _{2}}\widehat{\phi f_{\varepsilon }}\left( \xi -\eta \right)
\widehat{\phi g_{\varepsilon }}\left( \eta \right) d\eta +\int_{\Gamma
_{2}^{c}}\widehat{\phi f_{\varepsilon }}\left( \xi -\eta \right) \widehat{%
\phi g_{\varepsilon }}\left( \eta \right) d\eta \\
&=&I_{1}\left( \xi \right) +I_{2}\left( \xi \right)
\end{eqnarray*}
From (\ref{4.3}), $\exists C_{1}>0,\exists k_{1},k_{2}>0,\exists \varepsilon
_{1}>0$ such that $\forall \varepsilon \leq \varepsilon _{1},\forall \eta
\in \Gamma _{2},$
\begin{equation*}
\widehat{\phi f_{\varepsilon }}\left( \xi -\eta \right) \leq C_{1}\exp
\left( k_{1}\varepsilon ^{-\frac{1}{2s-1}}-k_{2}\left| \xi -\eta \right| ^{%
\frac{1}{s}}\right) ,
\end{equation*}
and from proposition \ref{ref3} $\exists C_{2}>0,\exists k_{3}>0,\forall
k_{4}>0,\exists \varepsilon _{2}>0,\forall \eta \in \mathbb{R}^{m},\forall
\varepsilon \leq \varepsilon _{2},$%
\begin{equation*}
\left| \widehat{\phi g_{\varepsilon }}\left( \eta \right) \right| \leq
C_{2}\exp \left( k_{2}\varepsilon ^{-\frac{1}{2s-1}}+k_{4}\left| \eta
\right| ^{\frac{1}{s}}\right)
\end{equation*}
Let $\gamma >0$ sufficiently small such that $\left| \xi -\eta \right| ^{%
\frac{1}{s}}\geq \gamma \left( \left| \xi \right| ^{\frac{1}{s}}+\left| \eta
\right| ^{\frac{1}{s}}\right) ,\forall \eta \in \Gamma _{2}.$ Hence for $%
\varepsilon \leq \min \left( \varepsilon _{1},\varepsilon _{2}\right) ,$%
\begin{equation*}
\left| I_{1}\left( \xi \right) \right| \leq C_{1}C_{2}\exp \left( \left(
k_{1}+k_{2}\right) \varepsilon ^{-\frac{1}{2s-1}}-k_{2}\gamma \left| \xi
\right| ^{\frac{1}{s}}\right) \int \exp \left( -k_{2}\gamma \left| \eta
\right| ^{\frac{1}{s}}+k_{4}\left| \eta \right| ^{\frac{1}{s}}\right) d\eta
\end{equation*}
take $k_{4}>k_{2}\gamma ,$ then
\begin{equation}
\left| I_{1}\left( \xi \right) \right| \leq C^{\prime }\exp \left(
k_{1}^{\prime }\varepsilon ^{-\frac{1}{2s-1}}-k_{2}^{\prime }\left| \xi
\right| ^{\frac{1}{s}}\right)  \label{4-4}
\end{equation}
Let $r>0,$
\begin{eqnarray*}
I_{2}\left( \xi \right) &=&\int_{\Gamma _{2}^{c}\cap \left\{ \left| \eta
\right| \leq r\left| \xi \right| \right\} }\widehat{\phi f_{\varepsilon }}%
\left( \xi -\eta \right) \widehat{\phi g_{\varepsilon }}\left( \eta \right)
d\eta +\int_{\Gamma _{2}^{c}\cap \left\{ \left| \eta \right| \geq r\left|
\xi \right| \right\} }\widehat{\phi f_{\varepsilon }}\left( \xi -\eta
\right) \widehat{\phi g_{\varepsilon }}\left( \eta \right) d\eta \\
&=&I_{21}\left( \xi \right) +I_{22}\left( \xi \right)
\end{eqnarray*}
Choose $r$ sufficiently small such that $\left\{ \left| \eta \right| ^{\frac{%
1}{s}}\leq r\left| \xi \right| ^{\frac{1}{s}}\right\} \Longrightarrow \xi
-\eta \notin \Gamma _{1}.$ Then $\left| \xi -\eta \right| ^{\frac{1}{s}}\geq
\left( 1-r\right) \left| \xi \right| ^{\frac{1}{s}}\geq \left( 1-2r\right)
\left| \xi \right| ^{\frac{1}{s}}+\left| \eta \right| ^{\frac{1}{s}}$ ,
consequently $\exists C_{3}>0,\exists \lambda _{1},\lambda _{2},\lambda
_{3}>0,\exists \varepsilon _{3}>0$ such that $\forall \varepsilon \leq
\varepsilon _{3},$
\begin{eqnarray*}
\left| I_{21}\left( \xi \right) \right| &\leq &C_{3}\exp \left( \lambda
_{1}\varepsilon ^{-\frac{1}{2s-1}}\right) \int \exp \left( -\lambda
_{2}\left| \xi -\eta \right| ^{\frac{1}{s}}-\lambda _{3}\left| \eta \right|
^{\frac{1}{s}}\right) d\eta \\
&\leq &C_{3}\exp \left( \lambda _{1}\varepsilon ^{-\frac{1}{2s-1}}-\lambda
_{2}^{\prime }\left| \xi \right| ^{\frac{1}{s}}\right) \int \exp \left(
-\lambda _{3}^{\prime }\left| \eta \right| ^{\frac{1}{s}}\right) d\eta \\
&\leq &C_{3}^{\prime }\exp \left( \lambda _{1}\varepsilon ^{-\frac{1}{2s-1}%
}-\lambda _{2}^{\prime }\left| \xi \right| ^{\frac{1}{s}}\right)
\end{eqnarray*}
If $\left| \eta \right| ^{\frac{1}{s}}\geq r\left| \xi \right| ^{\frac{1}{s}%
},$ we have $\left| \eta \right| ^{\frac{1}{s}}\geq \dfrac{\left| \eta
\right| ^{\frac{1}{s}}+r\left| \xi \right| ^{\frac{1}{s}}}{2}$ , and then $%
\exists C_{4}>0,\exists \mu _{1},\mu _{3}>0,\forall \mu _{2}>0,\exists
\varepsilon _{4}>0$ such that $\forall \varepsilon \leq \varepsilon _{4},$
\begin{eqnarray*}
\left| I_{21}\left( \xi \right) \right| &\leq &C_{4}\exp \left( \mu
_{1}\varepsilon ^{-\frac{1}{2s-1}}\right) \int \exp \left( \mu _{2}\left|
\xi -\eta \right| ^{\frac{1}{s}}-\mu _{3}\left| \eta \right| ^{\frac{1}{s}%
}\right) d\eta \\
&\leq &C_{4}\exp \left( \mu _{1}\varepsilon ^{-\frac{1}{2s-1}}\right) \int
\exp \left( \mu _{2}\left| \xi -\eta \right| ^{\frac{1}{s}}-\mu _{3}^{\prime
}\left| \eta \right| ^{\frac{1}{s}}-\mu _{3}^{\prime ^{\prime }}\left| \xi
\right| ^{\frac{1}{s}}\right) d\eta ,
\end{eqnarray*}
if take $\mu _{2}<\frac{\mu _{3}^{\prime }}{2}\left( 1+\frac{1}{r}\right) ,$
we obtain
\begin{equation*}
\left| I_{21}\left( \xi \right) \right| \leq C_{4}^{\prime }\exp \left(
k_{3}^{\prime }\varepsilon ^{-\frac{1}{2s-1}}-\mu _{3}^{\prime ^{\prime
}}\left| \xi \right| ^{\frac{1}{s}}\right) ,
\end{equation*}
which finish the proof
\end{proof}


\begin{thebibliography}{99}
\bibitem{A-R}  A. B. Antonevich, Ya. V. Radyno. On general method of
constructing algebras of new generalized functions. Soviet. Math. Dokl.,
vol. 43:3, ( 1991), 680-684

\bibitem{Ben-Bouz}  K. Benmeriem, C. Bouzar. Colombeau generalized functions
and solvability of differential operators. To appear in Zeitschrift f\"{u}
Analysis und ihre Anwendungen, (2006)

\bibitem{Biag}  H. A. Biagioni. A nonlinear theory of generalized functions.
Lecture Notes Math. 1421, Springer, 1990

\bibitem{Col1}  J. F. Colombeau. New generalized functions and
multiplication of distributions, North Holland, 1983

\bibitem{Col2}  J. F. Colombeau. Elementary introduction to new generalized
functions. North Holland, 1984

\bibitem{Col3}  J. F. Colombeau. Multiplication of Distributions: a tool in
mathematics numerical engineering and theorical physics. Lecture Notes in
Math. 1532, Springer-Verlag, 1992

\bibitem{DHMV}  A. Delcroix, M. Hasler, J.-A. Marti, V. Valmorin (Editors).
Nonlinear algebraic analysis and applications. Cambridge Scientific
Publishers. 2004

\bibitem{DHPV}  A. Delcroix, M. Hasler, S. Pilipovic, V. Valmorin.
Embeddings of ultradistributions and periodic hyperfunctions in Colombeau
type algebras through sequence space. Math. Proc. Cambr. Phil. Soc., 137:3,
(2004), 697-708

\bibitem{Eid}  A. Eida. On the microlocal decomposition of
ultradistributions and ultradifferentiable functions. Bull. Acad. Serbe Sc.
Arts, T. CXXII, (2001), 75-105

\bibitem{GKOS}  M. Grosser, M. Kuzinger, M. Oberguggenberger, R. Steinbauer.
Geometric theory of generalized functions with applications to general
relativity, Kluwer Publishing, 2001

\bibitem{Gel}  I. M. Guelfand, G. E. Shilov. Generalized functions, Vol. 2,
Academic Press, 1967

\bibitem{Gram}  T. Gramchev. Nonlinear maps in space of distributions, Math.
Z., 209, (1992), 101-114

\bibitem{GHKO}  M. Grosser, G. H\"{o}rmann, M. Kunzinger and M.
Oberguggenberger (Editors) Nonlinear Theory of Generalized Functions.
Chapman and Hall. 1999

\bibitem{Hor}  L. H\"{o}rmander. Distributions theory and Fourier analysis,
Springer-Verlag 1983

\bibitem{Hor-Kun}  G. H\"{o}rmann, M. Kunzinger. Microlocal properties of
basic operations in Colombeau algebras. J. Math. Anal. Appl., 261, (2001),
254-270

\bibitem{HorObPil}  G. H\"{o}rmann, M. Oberguggenberger, S. Pilipovi\'{c}.
Microlocal hypoellipticity of linear differential operators with generalized
functions as coefficients, Trans.\ Amer. Math. Soc., (2005)

\bibitem{Kom}  H. Komatsu. Ultradistributions I, J. Fac. Sci. Univ. Tokyo,
Sect. IA, 20, (1973), 25-105

\bibitem{LM}  J. L. Lions, E. Magenes. Non-Homogeneous Boundary Value
Problems and Applications , Vol.3, Springer-Verlag,1973

\bibitem{Mart}  J. A. Marti. ($\mathcal{C},\mathcal{E},\mathcal{P}$)-Sheaf
structures and applications. In M. Grosser, G. H\"{o}rmann, M. Kunzinger and
M. Oberguggenberger (Editors) Nonlinear Theory of Generalized Functions. pp
175-186. Chapman and Hall. 1999

\bibitem{NPS}  M. Nedeljkov, S. Pilipovic, D. Scarpal\'{e}zos. The linear
theory of Colombeau generalized functions, Longman Scientific \& Technica,
1998

\bibitem{Ober}  M. Oberguggenberger. Multiplication of distributions and
applications to partial differential equations, Longman Scientific \&
Technical, 1992

\bibitem{Pil}  S. Plilipovic. Microlocal properties of ultradistributions.
Composition and kernel type operators, Publ. Institut Math., T. 64, (1998),
85-97

\bibitem{Pilip}  S. Plilipovic, D. Scarpalezos. Colombeau generalized
ultradistributions. Math. Proc. Camb. Phil. Soc., 130, (2001), 541-553

\bibitem{Rod}  L. Rodino. Linear partial differential operators in Gevrey
spaces. World Scientific. 1993

\bibitem{Schw}  L. Schwartz. Sur l'impossibilit\'{e} de la multiplication
des distributions. C. R. Acad. Sci., 239, (1954), 847-848

\bibitem{Schw2}  L. Schwartz. Th\'{e}orie des distributions. Hermann, second
edition, 1966
\end{thebibliography}
\end{document}